\documentclass[runningheads]{llncs}
\usepackage{graphicx}
\usepackage{float}
\usepackage{mathtools}
\usepackage{amssymb, amsmath, latexsym}
\usepackage{nicefrac}
\usepackage{hyperref}
\usepackage{array}
\usepackage{hhline}
\usepackage{multirow}
\usepackage{wrapfig}
\usepackage{lipsum}
\usepackage{pifont}
\usepackage[colorinlistoftodos,bordercolor=orange,backgroundcolor=orange!20,linecolor=orange,textsize=scriptsize]{todonotes}
\usepackage{algorithm}\usepackage{algpseudocode}

\def \RR {\mathbb R}
\def \R {\mathbb R}
\def\eqdef{\overset{\text{def}}{=}}

\newcommand{\EndProof}{\begin{flushright}$\square$\end{flushright}}

\newcommand{\EE}{\mathbf{E}}

\def\R{\mathbb{R}}

\def\R{\mathbb R}

\def\EE{\mathbb E}

\def\la{\langle}
\def\ra{\rangle}

\newcolumntype{M}[1]{>{\centering\arraybackslash}m{#1}}
\newcolumntype{N}{@{}m{0pt}@{}}

\begin{document}
\title{Zeroth-Order Algorithms for Smooth Saddle-Point Problems
\thanks{The research of A. Sadiev, A. Beznosikov and P. Dvurechensky in Section 4  was supported by the Russian Science Foundation (project 21-71-30005). The research of A. Gasnikov in Section 5 was partially supported by RFBR, project number 18-29-03071 mk. This work was partially conducted while A. Sadiev and A. Beznosikov were on the project internship in Sirius University of Science and Technology.}}
\titlerunning{Zeroth-Order Algorithms for Smooth Saddle-Point Problems}
%
\author{ Abdurakhmon Sadiev\inst{1}\and
Aleksandr Beznosikov\inst{1,2}\and
Pavel Dvurechensky\inst{3,4} \and
Alexander Gasnikov\inst{1,4,5}}
\authorrunning{A. Sadiev, A. Beznosikov  et al.}
%
\institute{Moscow Institute of Physics and Technology, Russia \and
Higher School of Economics, Russia \and 
Weierstrass Institute for Applied Analysis and Stochastics, Germany \and
Institute for Information Transmission Problems RAS, Russia  \and
Caucasus Mathematical Center, Adyghe State University, Russia
}
\maketitle              
\begin{abstract}
Saddle-point problems have recently gained an increased attention from the machine learning community, mainly due to applications in training Generative Adversarial Networks using stochastic gradients. At the same time, in some applications only a zeroth-order oracle is available.
In this paper, we propose several algorithms to solve stochastic smooth (strongly) convex-concave saddle-point problems using zeroth-order oracles, and estimate their convergence rate and its dependence on the dimension $n$ of the variable. In particular, our analysis shows that in the case when the feasible set is a direct product of two simplices, our convergence rate for the stochastic term is only by a $\log n$ factor worse than for the first-order methods.
We also consider a mixed setup and develop 1/2th-order methods which use zeroth-order oracle for the minimization part and first-order oracle for the maximization part. 
 Finally, we demonstrate the practical performance of our zeroth-order and 1/2th-order methods on practical problems.


\end{abstract}

\section{Introduction}

Zeroth-order or derivative-free methods \cite{rosenbrock1960automatic,fabian1967stochastic,brent1973algorithms,spall2003introduction,conn2009introduction} are well known in optimization in application to problems with unavailable or computationally expensive gradients. In particular, the framework of derivative-free methods turned out to be very fruitful in application to different learning problems such as online learning in the bandit setup \cite{bubeck2012regret} and reinforcement learning \cite{salimans2017evolution,choromanski2018structured,fazel2018global}, which can be considered as a particular case of simulation optimization \cite{fu2015handbook,shashaani2018ASTRO}. 
We study stochastic derivative-free methods in a two-point feedback situation, meaning that two observations of the objective per iteration are available. This setting was considered for optimization problems by \cite{agarwal2010optimal,duchi2015optimal,shamir2017optimal} in the learning community and by \cite{nesterov2017random,stich2011linear,ghadimi2013stochastic,ghadimi2016mini-batch,gasnikov2016gradient-free,dvurechensky2020accelerated} in the optimization community.

In this paper we go beyond the setting of optimization problems and consider convex-concave saddle-point problems for which partial derivatives of the objective are not available, which forces to use derivative-free methods. Saddle-point problems are tightly connected with equilibrium \cite{facchinei2007finite} and game problems \cite{basar1998dynamic} in many applications, e.g., economics \cite{morgenstern1953theory}, with tractable reformulations of non-smooth optimization problems \cite{Nemirovski2004}, and with variational inequalities \cite{harker1990finite}. Gradient methods for saddle-point problems are an area of intensive study in the machine learning community in application to training of Generative Adversarial Networks \cite{goodfellow2014generative}, and other adversarial models \cite{madry2018towards}, as well as to robust reinforcement learning \cite{pinto2017robust}. In the latter two applications, gradients are often unavailable, which motivates the application of zeroth-order methods to the respective saddle-point problems. Moreover, this also motivates 1/2th-order methods, when the training of the network is made via stochastic gradient method with backpropagation, and adversarial examples, which are generated to force the network to give incorrect prediction, are generated by zeroth-order methods.    
Another application area for zeroth-order methods are
Adversarial Attacks \cite{goodfellow2014explaining,tramr2017ensemble}, in particular the Black-Box Adversarial Attacks \cite{conf/cvpr/NarodytskaK17}. The goal is not only to train the network, but to find also a perturbation of the data in such a way that the network outputs wrong prediction. Then the training is repeated to make the network robust to such adversarial examples. Since the attacking model does not have access to the architecture of the main network, but only to the input and output of the network, the only available oracle for the attacker is the zeroth-order oracle for the loss function. As it is shown in \cite{croce2018randomized,ye2018hessianaware,croce2019scaling}, this approach allows to obtain the same quality of robust training as the more laborious methods of Adversarial Attacks, but faster in up to a factor of three in terms of the training time \cite{Chen_2017}. 

Gradient methods for saddle-point problems are a well studied area with the classical algorithm being the extra-gradient method \cite{Korpelevich1976TheEM}. It was later generalized to the non-Euclidean geometry in the form of Mirror Descent \cite{nemirovski} and Mirror-Prox \cite{Nemirovski2004}. These methods are designed for a more general problem of solving variational inequalities. 
There are also direct methods for saddle-point problems such as gradient descent ascent \cite{nedic2009subgradient} or primal-dual hybrid gradient method \cite{chambolle2011first-order} for saddle-point problems with bilinear structure. On the contrary, the theory of zeroth-order methods for saddle-point problems seems to be underdeveloped in the literature. We give a more detailed overview of such methods and explain our contribution in comparison with the literature below.

\subsection{Our contribution and related works}

In the first part of the work, we present \textit{zeroth-order} variants of Mirror-Descent \cite{nemirovski} and  Mirror-Prox \cite{juditsky2008solving} methods for  \textit{stochastic saddle-point problems} in convex-concave and strongly convex-concave cases. We consider various concepts of zeroth-order oracles and various concepts of noise. Also we introduce a new class of smooth saddle-point problems -- firmly  smooth.

In the particular case of deterministic problems, our methods have a linear rate in the smooth strongly-convex-strongly-concave case, and sublinear rate $\mathcal{O}({1/N})$ in the convex-concave case, where $N$ is the number of iterations.
One can note that in some estimates, there is a factor of the problem's dimension $n$, but somewhere $n^{2/q}$. This factor $q$ depends on geometric setup of our problem and gives a benefit when we work in the Hölder, but non-Euclidean case (use non-Euclidean prox), i.e. $\| \cdot\| = \|\cdot\|_p$ and $p \in [1;2]$, then $\| \cdot \|_* = \| \cdot\|_q$, where $\nicefrac{1}{p} +\nicefrac{1}{q} = 1$. Then $q$ takes values from $2$ to $\infty$, in particular, in the Euclidean case $q=2$, but when the optimization set is a simplex, $q = \infty$. 
(see Table \ref{summary_1} for a comparison of the oracle complexity with zeroth-order methods for saddle-point problems in the literature and provided by our methods).
\renewcommand{\arraystretch}{1.5}
\renewcommand{\tabcolsep}{7pt} 
\begin{table}[h!]
\begin{center}
\begin{tabular}{cccc}
\hline
\textbf{Method} & \textbf{Assumptions} & \textbf{Complexity in deterministic setup}  \\\hline
{\tt ZO-GDMSA} \cite{wang2020zerothorder}  &  NC-SC, UCst-Cst, S  & $\mathcal{\tilde O} \left(\frac{n \kappa^2 }{\varepsilon^2}\right)$ \\\hline
{\tt ZO-Min-Max} \cite{liu2019minmax}  &  NC-SC, Cst-Cst, S  & $\mathcal{\tilde O} \left(\frac{n }{\varepsilon^6}\right)$ \\\hline
{\tt zoSPA} \cite{beznosikov_sadiev_gasnikov} & C-C, Cst-Cst, BG & $\mathcal{O}\left(n^{\nicefrac{2}{q}} \frac{M^2 D^2 }{\varepsilon^2}\right)$ \\\hline
 [\textbf{Alg 1 and 3}]  & SC-SC, Cst-Cst, S & $\mathcal{\tilde O}\left(\min\left[n^{\nicefrac{2}{q}} \kappa^2 , n \kappa  \right] \cdot \log\left(\frac{1}{\varepsilon}\right) \right)$ \\\hline
  [\textbf{Alg 2}] & C-C, Cst-Cst, S & $\mathcal{\tilde O}\left(n \frac{LD^2}{\varepsilon}\right)$ \\\hline
  [\textbf{Alg 1}] & C-C, Cst-Cst, FS & $\mathcal{\tilde O}\left(  n^{\nicefrac{2}{q}} \frac{L^2D^2}{\varepsilon} \right)^*$ \\\hline
\end{tabular} 
\end{center}
\caption{Comparison of oracle complexity in deterministic setup of different zeroth-order methods with different assumptions on target function $f(x,y)$: C-C -- convex-concave, SC-SC -- strongly-convex-strongly-concave, NC-SC -- nonconvex-strongly-concave; Cst -- optimizaation set is constrained, UCst -- unconstrained; S - smooth, FS - firmly smooth (see \eqref{SP_fc}), BG - bounded gradients. Here $\varepsilon$ means the accuracy of the solution, $D$ -- the diameter of the optimization set, $\mu$ -- strong convexity constant (see \eqref{SP_sm}), $L$ -- smoothness constant (see \eqref{SP_c}), $\kappa = \nicefrac{L}{\mu}$, $M$ -- bound of the gradient ($\| \nabla_x f(x,y)\|_2 \leq M$, $\| \nabla_y f(x,y)\|_2 \leq M$), $n$ -- the sum of the dimensions of the variables $x$ and $y$, $q = 2$ for the Euclidean case and $q = \infty$ for setup of $\|\cdot\|_1$-norm.\hspace{\textwidth}*convergence on $\frac{1}{N}\sum^N_{k = 1}\EE\left[\|F(x_k, y_k) - F(x^*, y^*)\|^2_2\right]$, where  $F(x,y) = (\nabla_x f(x,y),-\nabla_y f(x,y))$.}
\label{summary_1}
\end{table}

Our theoretical analysis shows that the zeroth-order methods has the same sublinear convergence rate in the stochastic part as the first-order method:  $\mathcal{O}({1/\sqrt{N}})$ in convex-concave case and $\mathcal{O}({1/N})$ in strongly-convex-strongly-concave case. (see Table \ref{summary_2} for a comparison of the oracle complexity in the stochastic part for first-order methods and available zeroth-order methods for stochastic saddle-point problems).
\renewcommand{\arraystretch}{1.5}
\renewcommand{\tabcolsep}{5pt} 
\begin{table}[h!]
\begin{center}
\begin{tabular}{cccc}
\hline
\textbf{Method} & \textbf{Order} & \textbf{Assumptions} & \textbf{Complexity for stochastic part}  \\\hline
{\tt EGMP} \cite{juditsky2008solving}  & 1st   & C-C, Cst-Cst, S  & $\mathcal{O}\left(\frac{\sigma^2 D^2}{\varepsilon^2}\right)$      \\\hline
{\tt PEG} \cite{hsieh2019convergence} & 1st  & SC-SC, Cst-Cst, S & $\mathcal{O}\left(\frac{\sigma^2}{\mu^2 \varepsilon}\right)$   \\\hline
{\tt ZO-SGDMSA}\cite{wang2020zerothorder} & 0th   & NC-SC, UCst-Cst, S & $\mathcal{\tilde O} \left(\frac{\kappa^2 n\sigma^2}{\varepsilon^4}\right)$  \\\hline
     [\textbf{Alg 1}]  & 0th   & SC-SC, Cst-Cst, S   &  $\mathcal{O}\left(\frac{n^{2/q}\sigma^2}{\mu^2 \varepsilon}\right)$  \\\hline
     [\textbf{Alg 2}]  & 0th   & C-C, Cst-Cst, S   &  $\mathcal{O}\left(\frac{n\sigma^2 D^2}{\varepsilon^2}\right)$  \\\hline
     [\textbf{Alg 1}]  & 0th   & C-C, Cst-Cst, FS   &  $\mathcal{O}\left(\frac{n^{2/q}\sigma^2 D^2}{\varepsilon^2}\right)$   \\\hline
\end{tabular} 
\end{center}
\caption{Comparison of oracle complexity for stochastic part of different first- and zeroth-order methods with different assumptions on $f(x,y)$: see 
notation in Table \ref{summary_1}. Here $\sigma^2$ -- the bound of variance (see \eqref{ran_dir1}).}
\label{summary_2}
\end{table}

The second part of the work is devoted to the use of a mixed order oracle, i.e. a zeroth-order oracle in one variable and a first-order oracle for the other. First, we analyze a special case when such an approach is appropriate - the Lagrange multiplier method. Then we also present a general approach for this setup. The idea of using such an oracle is found in the in literature \cite{aleks2019derivativefree}, but for the composite optimization problem.

As mentioned above, all theoretical results are tested in practice on a classical bilinear problem.


\section{Problem setup and assumptions}\label{sec:statement}

We consider a saddle-point problem: 
\begin{equation}
    \label{SP}
    \min_{x \in \mathcal{X}} \max_{y \in \mathcal{Y}} f(x,y),
\end{equation}
where $\mathcal{X} \subset \R^{n_x}$ and $\mathcal{Y} \subset \R^{n_y}$ are convex compact sets. For simplicity, we introduce the set $\mathcal{Z} = \mathcal{X} \times \mathcal{Y}$, $z = (x,y)$ and the operator $F$:
\begin{equation}
\label{opSP}
    F(z) = F(x,y) = \begin{pmatrix}
\nabla_x f(x,y)\\
-\nabla_y f(x,y)
\end{pmatrix}.
\end{equation}
We focus on the case when we do not have access to the values of $\nabla_x f(x,y)$ and $\nabla_y f(x,y)$, but we have access to the inexact zeroth-order oracle, i.e. inexact values of the objective $f(x,y)$. The inexactness in the zeroth-order oracle includes  stochastic noise and unknown bounded noise, which can be of an adversarial nature. More precisely, we have access to the values $\tilde f(z, \xi)$ such that $\tilde f(z, \xi) = f(z, \xi) + \delta(z)$ and
\begin{eqnarray}
\label{ran_dir1}
    \EE[f(z,\xi)] = f(z), ~~\EE[F(z,\xi)] = F(z), \nonumber\\
    \EE[\|F(z,\xi) - F(z)\|^2_2] \leq \sigma^2, ~~ |\delta(z)| \leq \Delta.
\end{eqnarray}

We consider two types of approximations for $F(z)$ based on the available observations of $\tilde f(z, \xi)$.

\textbf{Random direction oracle.} In this strategy, the vectors $e_x, e_y$ are generated uniformly on the unit Euclidean sphere, i.e. $e_x\in \mathcal{RS}^2_{n_x}(1)$ and $e_y\in \mathcal{RS}^2_{n_y}(1)$. And 
\begin{eqnarray}
    \label{ran_dir}
    g_d(z,e,\tau, \xi) = \frac{n}{\tau} \left(
    \begin{array}{c}
    \left(\tilde f(x + \tau e_x, y,  \xi) -  \tilde f(x, y, \xi)  \right)e_x\\
    \left(\tilde f(x, y, \xi) - \tilde f(x, y + \tau e_y,  \xi)   \right)e_y \\
    \end{array}
    \right),
\end{eqnarray}
where $\tau >0$ is called smoothed parameter and $n = n_x + n_y + 1$.

\textbf{Full coordinates oracle.} Here we consider a standard orthonormal basis $\{h_1, \ldots, h_{n_x + n_y}\}$ and construct an approximation for the operator $F$ in the following form:
\begin{eqnarray}
    \label{full_coor}
    g_f(z,h,\tau, \xi) &=& \frac{1}{{\tau}}\sum\limits_{i=1}^{n_x}\left(\tilde f(z+\tau h_i, \xi) - \tilde f(z, \xi) \right)h_i \nonumber\\ &&+
    \frac{1}{{\tau}}\sum\limits_{i=n_x + 1}^{n_x + n_y}\left(\tilde f(z, \xi)  - \tilde f(z+\tau h_i, \xi) \right)h_i.
\end{eqnarray}
In this concept, we need to call $\tilde f$ oracle $n_x + n_y + 1$ times, whereas in the previous case only 3 times.


\section{Notation and Definitions}\label{sec:notation}

We use $\la x,y \ra \eqdef \sum_{i=1}^nx_i y_i$ to define inner product of $x,y\in\R^n$ where $x_i$ is the $i$-th component of $x$ in the standard basis in $\R^n$. Hence we get the definition of $\ell_2$-norm in $\R^n$ in the following way $\|x\|_2 \eqdef \sqrt{\la x, x \ra}$. We define $\ell_p$-norms as $\|x\|_p \eqdef \left(\sum_{i=1}^n|x_i|^p\right)^{\nicefrac{1}{p}}$ for $p\in(1,\infty)$ and for $p = \infty$ we use $\|x\|_\infty \eqdef \max_{1\le i\le n}|x_i|$. The dual norm $\|\cdot\|_q$ for the norm $\|\cdot\|_p$ is defined in the following way: $\|y\|_q \eqdef \max\left\{\la x, y \ra\mid \|x\|_p \le 1\right\}$. Operator $\EE[\cdot]$ is full mathematical expectation and operator $\EE_\xi[\cdot]$ express conditional mathematical expectation. 

As stated above, during the course of the paper we will work in an arbitrary norm $\|\cdot\| = \|\cdot\|_p$, where $p \in [1;2]$. And its conjugate $\|\cdot \|_* = \|\cdot\|_q$ with $q \in [2; +\infty)$ and $\nicefrac{1}{p} + \nicefrac{1}{q} = 1$. Some assumptions will be made later in the Euclidean norm - we will write this explicitly $\|\cdot\|_2$. 

\textbf{Definition 1.}
Function $d(z): \mathcal{Z} \to \mathbb{R}$ is called prox-function if $d(z)$ is $1$-strongly convex w.r.t. $\| \cdot \|$-norm and differentiable on $\mathcal{Z}$ function. 

\textbf{Definition 2.}
Let $d(z): \mathcal{Z} \to \mathbb{R}$ is prox-function. For  any  two  points  $z,w \in \mathcal{Z}$ we define Bregman divergence $V_z(w)$ associated with $d(z)$ as follows: 
\begin{equation*}
    \label{Bregman}
    V_z(w) = d(z) - d(w) - \langle \nabla d(w), z - w \rangle.
\end{equation*}


\textbf{Definition 3.}
Let $V_z(w)$ Bregman divergence. For all $x \in \mathcal{Z}$ define prox-operator of $\xi$:
\begin{equation*}
    \text{prox}_x (\xi) = \text{arg}\min_{y \in \mathcal{Z}} \left(V_x(y) + \langle \xi , y \rangle \right).
\end{equation*}

Next we present the assumptions that we will use in the convergence analysis.

\textbf{Assumption 1.} The set $\mathcal{Z}$ is bounded w.r.t $\|\cdot\|$ by constant $D_p$, i.e.
\begin{equation}
    \label{b}
    V_{z_1}(z_2) \leq D^2_p,~~ \forall z_1, z_2 \in \mathcal{Z}.
\end{equation}


\textbf{Assumption 2.} $f(x,y)$ is convex-concave. It means that $f(\cdot, y)$ is convex for all $y$ and  $f(x, \cdot)$ is concave for all $x$.

\textbf{Assumption 2(s).} $f(x,y)$ is strongly-convex-strongly-concave. It means that $f(\cdot, y)$ is strongly-convex for all $y$ and  $f(x, \cdot)$ is strongly-concave for all $x$ w.r.t. $V_{\cdot}(\cdot)$, i.e. for all $x_1, x_2 \in \mathcal{X}$ and for all $y_1, y_2 \in \mathcal{Y}$ we have
\begin{eqnarray}
\label{SP_sm}
f(x_1,y_2) &\geq& f(x_2, y_2) + \langle\nabla_x f(x_2, y_2) , x_1 - x_2\rangle \nonumber\\ &&+ \frac{\mu}{2}\left(V_{(x_2, y_2)}(x_1, y_2) + V_{(x_1, y_2)}(x_2, y_2)\right), \nonumber\\
-f(x_2,y_1) &\geq& -f(x_2, y_2) + \langle-\nabla_y f(x_2, y_2) , y_1 - y_2\rangle \nonumber\\ &&+ \frac{\mu}{2}\left(V_{(x_2, y_2)}(x_2, y_1) + V_{(x_1, y_1)}(x_2, y_2)\right) .
\end{eqnarray}

\textbf{Assumption 3.} $f(x,y, \xi)$ is $L(\xi)$-Lipschitz continuous w.r.t $\|\cdot\|_2$, i.e. for all $x_1, x_2 \in \mathcal{X}, y_1, y_2 \in \mathcal{Y}$ and $\xi$
\begin{eqnarray}
\label{SP_c}
\left\|\begin{pmatrix}
\nabla_x f(x_1,y_1, \xi)\\
-\nabla_y f(x_1,y_1, \xi)
\end{pmatrix} -\begin{pmatrix}
\nabla_x f(x_2,y_2, \xi)\\
-\nabla_y f(x_2,y_2, \xi)
\end{pmatrix} \right\|_2  \leq L(\xi) \left\| \begin{pmatrix}
x_1\\
y_1
\end{pmatrix} - \begin{pmatrix}
x_2\\
y_2
\end{pmatrix} \right\|_2.\end{eqnarray}

\textbf{Assumption 3(f).} $f(x,y)$ s $L$-firmly Lipschitz continuous w.r.t $\|\cdot\|_2$, i.e.
for all $x_1, x_2 \in \mathcal{X}, y_1, y_2 \in \mathcal{Y}$
\begin{eqnarray}
\label{SP_fc}
\left\|\begin{pmatrix}
\nabla_x f(x_1,y_1, \xi)\\
-\nabla_y f(x_1,y_1, \xi)
\end{pmatrix} -\begin{pmatrix}
\nabla_x f(x_2,y_2, \xi)\\
-\nabla_y f(x_2,y_2, \xi)
\end{pmatrix} \right\|^2_2  & \nonumber\\
&\hspace{-7cm} \leq L(\xi) \left\langle \begin{pmatrix}
\nabla_x f(x_1,y_1, \xi)\\
-\nabla_y f(x_1,y_1, \xi)
\end{pmatrix} - \begin{pmatrix}
\nabla_x f(x_2,y_2, \xi)\\
-\nabla_y f(x_2,y_2, \xi)
\end{pmatrix}
,\begin{pmatrix}
x_1\\
y_1
\end{pmatrix} - \begin{pmatrix}
x_2\\
y_2
\end{pmatrix} \right\rangle.\end{eqnarray}

For \eqref{SP_c} and \eqref{SP_fc} we assume that exists $L_2$ such that $\EE[L^2(\xi)] \leq L^2_2$. For deterministic case $L_2 $ is equal to deterministic constant $L$ (without $\xi$).

By Cauchy-Schwarz, \eqref{SP_c} follows from \eqref{SP_fc}.
It is easy to see that the assumptions 4 and 4(f) above can be easily rewritten in a more compact form using $F(z)$. For assumption 3(s) it is more complicated:
\begin{lemma} \label{lemma2}
If $f(x,y)$ is $\mu$-strongly convex on $x$ and $\mu$-strongly concave on $y$ w.r.t $V_{\cdot}(\cdot)$, then for $F(z)$ we have
\begin{equation*}
    \langle F(z_1) - F(z_2), z_1 - z_2 \rangle \geq \frac{\mu}{2}\left( V_{z_1}(z_2) + V_{z_2}(z_1)\right),~~ \forall z_1, z_2 \in \mathcal{Z}.
\end{equation*}
\end{lemma}

Hereinafter, we do not present the proofs of lemmas and theorems in the main part of the paper -- see the corresponding parts of the appendix. And we can present some properties of oracles \eqref{ran_dir}, \eqref{full_coor}: 

\begin{lemma} \label{lemma1}
Let $e \in \mathcal{RS}^2(1)$, i.e. uniformly distributed on the unit Euclidean sphere. Randomness comes from independent variables $e$, $\xi$ and a point $z$. Norm $\|\cdot\|_* = \|\cdot\|_q$ satisfies $q \in [2;+\infty)$. We introduce the constant $\rho_n$:
\begin{equation*}
    \rho_n = \min\{q-1, 16 \log(n) -8 \}.
\end{equation*}
Then under Assumption 3 or 3(f) the following statements hold:
\begin{itemize}
   \item for Random direction oracle
    \begin{eqnarray*}
        \EE\left[\|g_d(z,e,\tau, \xi)\|^2_q\right] &\leq& 48 n^{2/q} \rho_n \EE\left[\| F(z) -F(z^*) \|^2_2 \right]+ 48 n^{2/q} \rho_n \| F(z^*) \|^2_2 \nonumber\\ 
    && + 48 n^{2/q} \rho_n \sigma^2 + 8 n^{2/q + 1}\rho_n L^2\tau^2  \nonumber \\
    && + 16 \frac{n^{2/q + 1}\rho_n \Delta^2}{\tau^2}, \\
        \left\| \EE [g_d(z,e,\tau, \xi)] - F(z)\right\|_q &\leq& 2 n^{1/q + 1/2} \sqrt{\rho_n} L\tau + 4 n^{1/q + 1/2} \sqrt{\rho_n} \frac{\Delta}{\tau};
    \end{eqnarray*}
    \item for Full coordinates oracle 
    \begin{eqnarray*}
        \EE\left[\|g_f(z,\tau, \xi) - F(z)\|^2_q\right] &\leq& 3\sigma^2 + 3nL_2^2 \tau^2 + \frac{6n\Delta^2}{\tau^2}, \\
        \|\EE\left[g_f(z,\tau, \xi)] - F(z)\right\|_q &\leq& \sqrt{n} L \tau + \frac{2\sqrt{n}\Delta}{\tau}.
    \end{eqnarray*}
\end{itemize}
\end{lemma}

\section{Zeroth-Order Methods}\label{sec:main_res1}
In this part, we present methods for solving problem \eqref{SP}, which use only the zeroth-order oracle. First of all, we want to consider the classic version of the Mirror-Descent algorithm. For theoretical and practical analysis of this\\
\begin{minipage}{0.45\textwidth}
\begin{algorithm}[H]
\caption{{\tt zoVIA}}
	\label{alg1}
\begin{algorithmic}
\State 
\noindent {\bf Input:} $z_0$, $N$, $\gamma$, $\tau$.
\State Choose $\text{grad}$ to be either $g_d$ or $g_f$.
\For {$k=0,1, 2, \ldots, N$ }
    \State Sample indep. $e_k$, $\xi_k$.
    \State $d_{k} =  \text{grad}(z_{k}, e_{k}, \tau, \xi_k)$.
    \State $z_{k+1} = \text{prox}_{z_k}(\gamma \cdot d_{k})$.
\EndFor
\State 
\noindent {\bf Output:} $z_{N+1}$ or $\bar z_{N+1}$.
\end{algorithmic}
\end{algorithm}
\end{minipage}
\begin{minipage}{0.05\textwidth}
\end{minipage}
\begin{minipage}{0.53\textwidth}
algorithm in the non-smooth case, but with a bounded gradient, see \cite{nemirovski}(first order), \cite{beznosikov_sadiev_gasnikov}(zero order). The main problem of this approach is that it is difficult to analyze in the case when $f$ is convex-concave and Lipschitz continuous (Assumptions 2 and 3). But in practice, this algorithm does not differ much from its counterparts, which will be given below. Let us analyze this algorithm in convex-concave and strongly-convex-strongly-concave cases with Random direction oracle:
\end{minipage}

\begin{theorem} By Algorithm 1 with Random direction oracle
\begin{itemize}
    \item under Assumptions 1, 2, 3(f) and with $\gamma \leq \frac{1}{48 n^{\nicefrac{2}{q}}\rho_nL}$, we get
    \begin{eqnarray*}
        \frac{1}{N}\sum^N_{k = 1}\EE\left[\|F(z_k) - F(z^*)\|^2_2\right] &\leq&\frac{2LD^2_p}{\gamma N} + 48 \gamma n^{2/q} \rho_nL\left(\|F(z^*)\|^2_2 + \sigma^2\right) \nonumber\\&& + 8\gamma n^{2/q + 1} \rho_n L\left(   L_2^2 \tau^2 + 2 \frac{\Delta^2}{\tau^2}\right) \nonumber\\&& + 8 n^{1/q + 1/2} \sqrt{\rho_n}LD_p\left(L\tau + \frac{2\Delta}{\tau}\right);
    \end{eqnarray*}
    \item under Assumptions 1, 2(s), 3 and with $\gamma \leq \frac{\mu}{96 n^{\nicefrac{2}{q}}\rho_nL^2}$:
    \begin{eqnarray*}
        \EE\left[V_{z_{N+1}}(z^*)\right]&\leq& V_{z_0}(z^*) \exp\left(-\frac{\mu^2 N}{400n^{\nicefrac{2}{q}}\rho_n L^2}\right) + \nonumber\\&&+ \frac{24n^{2/q} \rho_n}{\mu^2 N}\left(\|F(z^*)\|^2_2 + \sigma^2 \right)  \\&&+ \frac{4n^{2/q + 1} \rho_n }{\mu^2 N}\left(L_2^2 \tau^2 + 2 \frac{ \Delta^2}{\tau^2}\right)\\&& + \frac{4 n^{1/q + 1/2} \sqrt{\rho_n}  D_p}{\gamma \mu^2 N}\left(L\tau + \frac{2\Delta}{\tau}\right).
    \end{eqnarray*}
\end{itemize}
\end{theorem}

\textbf{Remark.} In the first statement of the Theorem, we used an unusual convergence criterion, it can be interpreted as follows: let as the output $\tilde z_{N}$ of the algorithm we choose a random point from $z_0$ to $z_N$. Then
$$\EE\left[\|F(\tilde z_{N})\|^2_2\right] = \frac{1}{N+1}\sum^N_{k = 0}\EE\left[\|F(z_k)\|^2_2\right].$$

In this theorem and below, we draw attention to the fact that in the main part of the convergence there is a deterministic constant $L$, and in the parts that are responsible for noise --  $L_2$ (see \eqref{SP_c},\eqref{SP_fc}).


\begin{corollary} For Algorithm 1
\begin{itemize}
    \item under Assumptions 1, 2, 3(f) and with $\gamma =\min\left\{ \frac{1}{48 n^{\nicefrac{2}{q}}\rho_nL}, \frac{D_p}{n^{\nicefrac{1}{q}}\sqrt{\rho_n}\sigma \sqrt{N}}\right\},$
    \begin{eqnarray*}
    \tau = \Theta \left( \min\left\{\frac{\varepsilon}{n^{1/q + 1/2} \sqrt{\rho_n}L^2D_p}, \max\left[ \sqrt{\frac{\varepsilon}{n L_2^2}}, \frac{\sigma}{\sqrt{n}L_2} \right]\right\}\right), \quad \Delta = \mathcal{O} \left(L_2 \tau^2\right),
    \end{eqnarray*}
    the  oracle complexity  (coincides with the number of iterations) to find $\varepsilon$-solution (in terms of the convergence criterion from Theorem 1) is
    \begin{eqnarray*}
    N = \mathcal{O}\left( \max\left\{\frac{n^{\nicefrac{2}{q}}\rho_n L^2 D^2_p }{\varepsilon}, \frac{n^{\nicefrac{2}{q}}\rho_n\sigma^2D^2_p}{ \varepsilon^2}\right\}\right).
    \end{eqnarray*}

    \item under Assumptions 1, 2(s), 3 and with $\gamma = \frac{\mu}{96 n^{\nicefrac{2}{q}}\rho_nL^2}$,
    \begin{eqnarray*}
    \tau = \Theta \left( \min\left\{ \max\left[\frac{\sqrt{\varepsilon} L}{L_2}, \frac{\sigma}{\sqrt{n}L_2} \right], \max\left[\frac{\varepsilon \mu}{n^{\nicefrac{1}{q} + \nicefrac{1}{2}} \sqrt{\rho_n} L D_p},\frac{\sigma^2 \mu}{n^{\nicefrac{1}{q} + \nicefrac{1}{2}} \sqrt{\rho_n} L^3 D_p}\right]\right\}\right),
    \end{eqnarray*}
    $\Delta = \mathcal{O} \left(L_2 \tau^2\right)$, the  oracle complexity  (coincides with the number of iterations) to find $\varepsilon$-solution (in terms of the convergence criterion from Theorem 1) can be bounded by
    \begin{eqnarray*}
    N = \widetilde{\mathcal{O}} \left( \max\left\{\frac{n^{\nicefrac{2}{q}}\rho_n L^2 }{\mu^2}\log\left(\frac{1}{\varepsilon}\right), \frac{n^{\nicefrac{2}{q}}\rho_n \sigma^2}{\mu^2 \varepsilon}\right\}\right).
\end{eqnarray*}

\end{itemize}

\end{corollary}

\textbf{Remark.} 
We analyze only Random direction oracle. The estimate of the oracle complexity with Full coordinate oracle has the same form with $q = 2$.

\vspace{\baselineskip}
Next, we consider a standard algorithm for working with smooth saddle-point problem. It builds on the extra-gradient method \cite{Korpelevich1976TheEM}. The idea of using this approach for saddle-point problems is not new \cite{juditsky2008solving}. It has both heuristic advantages (we forestall the properties of the gradient) as well as purely mathematical ones (a more clear theoretical analysis). We use two versions of this approach: classic and single call version from \cite{hsieh2019convergence}.\\
\begin{minipage}{0.53\linewidth}
     \begin{algorithm}[H]
	\caption{{\tt zoESVIA}}
	\label{alg2}
\begin{algorithmic}
\State 
\noindent {\bf Input:} $z_0$, $N$, $\gamma$, $\tau$.
\State Choose oracle $\text{grad}$ from $g_d, g_f$.
\For {$k=0,1, 2, \ldots, N$ }
    \State Sample indep. $e_k$, $e_{k+1/2}$, $\xi_k$, $\xi_{k+1/2}$.
    \State $d_{k} =  \text{grad}(z_{k}, e_{k}, \tau, \xi_k)$.
    \State $z_{k+1/2} = \text{prox}_{z_k}(\gamma \cdot d_{k})$.
    \State $d_{k+1/2} =  \text{grad}(z_{k+1/2}, e_{k+1/2}, \tau, \xi_{k+1/2})$.
    \State $z_{k+1} = \text{prox}_{z_k}(\gamma \cdot d_{k+1/2})$.
\EndFor
\State 
\noindent {\bf Output:} $z_{N+1}$ or $\bar z_{N+1}$.
\end{algorithmic}
\end{algorithm}
\end{minipage}
\begin{minipage}{0.01\linewidth}
\end{minipage}
\begin{minipage}{0.46\linewidth}
     \begin{algorithm}[H]
	\caption{{\tt zoscESVIA}}
	\label{alg3}
\begin{algorithmic}
\State 
\noindent {\bf Input:} $z_0$, $N$, $\gamma$, $\tau$.
\State Choose oracle $\text{grad}$ from $g_d, g_f$.
\For {$k=0,1, 2, \ldots, N$ }
    \State Sample independent $e_k$, $\xi_k$.
    \State Take $d_{k-1}$ from previous step.
    \State $z_{k+1/2} = \text{prox}_{z_k}(\gamma \cdot d_{k-1})$.
    \State $d_{k} =  \text{grad}(z_{k+1/2}, e_{k+1/2}, \tau, \xi_k)$.
    \State $z_{k+1} = \text{prox}_{z_k}(\gamma \cdot d_k)$.
\EndFor
\State 
\noindent {\bf Output:} $z_{N+1}$ or $\bar z_{N+1}$.
\end{algorithmic}
\end{algorithm}
\end{minipage}
Here $\bar z_{N+1} = \frac{1}{N+1} \sum_{i=0}^{N} z_{i+1/2}$.

Next, we will deal with the theoretical analysis of convergence:
\begin{theorem}\label{t2}
\begin{itemize}
    \item By Algorithm 2  with Full coordinates oracle under Assumptions 1, 2, 3 and with $\gamma \leq \nicefrac{1}{2L}$, we have
    \begin{eqnarray*}
        \EE\left[\varepsilon_{sad}(\bar z_{N+1}) \right] &\leq& \frac{2D_p^2}{\gamma N} + 11\gamma\left(n L_2^2 \tau^2 + \sigma^2 + 2 \frac{n\Delta^2}{\tau^2}\right) \nonumber\\
    &&+  2D_p\left( \sqrt{n} L \tau + \frac{2\sqrt{n}\Delta}{\tau}\right),
    \end{eqnarray*}
    where 
\begin{equation*}
    \varepsilon_{sad}(\bar z_{N+1}) = \max_{y' \in \mathcal{Y}} f(\bar x_{N+1}, y') - \min_{x' \in \mathcal{X}} f(x', \bar y_{N+1}), 
\end{equation*}
$\bar x_{N+1}$, $\bar y_{N+1}$ are defined the same way as $\bar z_{N+1}$.
    \item By Algorithm 3 with Full coordinates oracle under Assumptions 1, 2(s), 3 and with $p=2$ ($V_x(y) = \nicefrac{1}{2}\|x-y\|^2_2$), $\gamma \leq \nicefrac{1}{6L}$:
    \begin{eqnarray*}
        \EE\big[\|z_{N+1} - z^*\|^2_2\big] &\leq& \exp\left(- \frac{\mu N}{12 L}\right) \left(\|z_0 - z^*\|^2_2 + \|g_f(z_{0},\tau, \xi_{0}) - g_f(z_{0},\tau, \xi_{0})\|^2_2\right)\nonumber\\
    &&+ \frac{1}{\mu^2 N} 12\left(\sigma^2 + nL^2_2\tau^2+ \frac{2n\Delta^2}{\tau^2}\right) \nonumber\\
    &&+\frac{1}{\mu^2 N}\frac{4D_2}{\gamma} \left(\sqrt{n}L\tau+\frac{2\sqrt{n}\Delta}{\tau}\right).
    \end{eqnarray*}
\end{itemize}
\end{theorem}


\begin{corollary} Let $\varepsilon$ -- accuracy of the solution (in terms of the convergence criterion from Theorem 2). 
\begin{itemize}
    \item For Algorithm 2  with Full coordinates oracle under Assumptions 1, 2, 3 with $\gamma = \min\left\{\nicefrac{1}{2L}, \nicefrac{D_p}{(\sigma \sqrt{N})}\right\}$ and additionally 
    \begin{eqnarray*}
        \tau = \mathcal{O} \left( \min\left\{\frac{\varepsilon}{\sqrt{n}LD_2}, \max\left[\sqrt{\frac{\varepsilon L}{nL^2_2}}, \frac{\sigma}{\sqrt{n}L_2}\right]\right\}\right), \quad\Delta = \mathcal{O} \left(L_2 \tau^2\right),
    \end{eqnarray*}
    we have the number of iterations to find $\varepsilon$-solution
    \begin{eqnarray*}
    N = \mathcal{O}\left(\max\left\{ \frac{LD^2_2 }{\varepsilon}, \frac{\sigma^2 D^2_p }{\varepsilon^2}\right\}\right).
\end{eqnarray*}

\item For Algorithm 3 with Full coordinates oracle under Assumptions 1, 2(s), 3,  with $p=2$ ($V_x(y) = \nicefrac{1}{2}\|x-y\|^2_2$), $\gamma = \nicefrac{1}{6L}$ and additionally
\begin{eqnarray*}
    \tau = \mathcal{O} \left( \min\left\{\max\left[\sqrt{\frac{\varepsilon \mu L}{L_2^2}}, \frac{\sigma}{\sqrt{n} L_2}\right], \max\left[\frac{\mu \varepsilon}{\sqrt{n}L D_2}, \frac{\sigma^2}{\sqrt{n}L^2 D_2}\right] \right\}\right), 
\end{eqnarray*}
$\Delta = \mathcal{O} \left(L_2 \tau^2\right)$, the number of iterations to find $\varepsilon$-solution:
\begin{eqnarray*}
    N = \mathcal{\widetilde{O}}\left( \max\left\{\frac{L}{\mu}\log\left(\frac{1}{\varepsilon}\right), \frac{\sigma^2}{ \mu^2\varepsilon}\right\}\right).
\end{eqnarray*}

\end{itemize}

\end{corollary}

\textbf{Remark.} The oracle complexity for the Full coordinate oracle  is $n$ times greater than the number of iterations.

The analysis is carried out only for the Full coordinate oracle. The main problem of using Random Direction is that their variance is tied to the norm of the gradient; therefore, using an extra step does not give any advantages over Algorithm \ref{alg1}. A possible way out of this situation is to use the same direction $e$ within one iteration of Algorithm \ref{alg2} -- this idea is implemented in Appendix \ref{same direction} and in Practice part.
It is interesting how it work in practice, because in the non-smooth case \cite{beznosikov_sadiev_gasnikov} the gain by the factor $n^{2/q}$ can be obtained.

\section{$\nicefrac{1}{2}$-Order Methods}\label{sec:main_res}

In this section, we have access to a first-order oracle in one of the variables, and in the other -- only a zeroth-order oracle. For such a case, we suggest using an oracle of the form: 

\begin{equation*}
   \widetilde{g}(z, \tau) = \begin{pmatrix}
[grad(x, y)]_x\\
-\nabla_y f(x,y)
\end{pmatrix},
\end{equation*}
where $[grad(x, y)]_x$ -- one of the zeroth-order approximations on variable $x$: \eqref{ran_dir} or \eqref{full_coor}. Before proving the general case, we consider one illustrative example:

\subsection{Lagrange multiplier method}

Let $\mathcal{X} \subset \mathbb{R}^n$ be a convex, compact set and functions $f(x), g_1(x),\dots, g_m(x)$ be convex, smooth. We solve the following optimization problem:
\begin{eqnarray*}
    &&\min_{x \in \mathcal{X}} f(x) ,\nonumber\\ &&\text{s.t. } g_i(x)\leq0~ \forall i\in 1,\ldots m.
\end{eqnarray*}
A dual problem to the original one:
\begin{eqnarray*}
    \max_{\lambda \in \bot_m}\min_{x \in \mathcal{X}} \mathcal{L}(x, \lambda) = f(x)+ \langle\lambda, g(x)\rangle,
\end{eqnarray*}
where $\bot_m = \{y\in \mathbb{R}^m~|~y_i \geq 0\}$ -- a positive orthant, $\mathcal{L}(x, \lambda)$ -- a Lagrange function, $\lambda$ -- a Lagrange multiplier, $g(x) = (g_1(x),\dots, g_m(x)))^T$.
We got a saddle-point problem that we want to solve using the zeroth-order method, i.e. only function values are available. But it turns out that we have access to $\nabla_{\lambda}\mathcal{L}(x, \lambda) = g(x)$ completely free: when we build the "gradient"{} on $x$ using finite differences, we call the value for $g (x)$ and immediately get the gradient $\lambda$.

For such a problem, the oracle of the zero and first orders can be called the same number of times. In general, it is unprofitable to calculate the gradient as many times as the zeroth-order oracles and a slightly different result is obtained:


\subsection{Universal approach with Full gradient method}

Define Mixed oracle:
\begin{equation*}
   \widetilde{g}_f(z, \tau) = \begin{pmatrix}
[g_f(x, y)]_x\\
-\nabla_y f(x,y)
\end{pmatrix},
\end{equation*}
then
\begin{theorem} By Algorithm \ref{alg2} under assumption 1, 2, 3 with Mixed oracle $\widetilde{g}_f$ and $\gamma \leq \nicefrac{1}{2L}$, we get
\begin{eqnarray*}
     \EE\left[\varepsilon_{sad}(\bar z_N) \right] &\leq& \frac{D_p^2}{\gamma N} + 2 D_p\left(\sqrt{n_x} L_2 \tau + \frac{2\sqrt{n_x}\Delta}{\tau} \right)\nonumber\\ &&+ 9\gamma\left(\sigma^2 + n_x L_2^2 \tau^2 + \frac{2n_x\Delta^2}{\tau^2}\right) 
\end{eqnarray*}

\end{theorem}

\begin{corollary}
To get accuracy $\varepsilon$ (in terms of the convergence criterion from Theorem 2) in Algorithm 2 with Mixed oracle, under Assumptions 1, 2, 3, with $\gamma = \min\left\{\nicefrac{1}{2L}, \nicefrac{D_p}{(\sigma \sqrt{N})}\right\}$,
    \begin{eqnarray*}
        \tau = \mathcal{O} \left( \min\left\{\frac{\varepsilon}{\sqrt{n}LD_p}, \max\left[\sqrt{\frac{\varepsilon L}{nL^2_2}}, \frac{\sigma}{\sqrt{n}L_2}\right]\right\}\right), \quad\Delta = \mathcal{O} \left(L_2 \tau^2\right),
    \end{eqnarray*}
we need to call Full coordinates oracle for $x$
    \begin{eqnarray*}
    N = \mathcal{O}\left(\max\left\{ \frac{LD^2_p }{\varepsilon}, \frac{\sigma^2 D^2_p }{\varepsilon^2}\right\}\right)~~\text{times}.
    \end{eqnarray*}
\end{corollary}

\section{Practice part}

The main goal of our experiments is to compare the Algorithms 1,2,3 and 4 (see Appendix \ref{same direction}) described in this paper with Full coordinate and Random direction oracles. We consider the classical 
bilinear saddle-point problem on a probability simplex:
\begin{eqnarray}
\label{exp_pr_4}
    \min_{x\in \Delta_n}\max_{y\in \Delta_k} \left[  y^T Cx\right],
\end{eqnarray}
This problem is often referred to as a matrix game (see Part 5 in \cite{nemirovski}). Two players $X$ and $Y$ are playing. The goal of player $Y$ is to win as much as possible by correctly choosing an action from 1 to $k$, the goal of player $X$ is to minimize the gain of player $X$ using his actions from $1$ to $n$. Each element of the matrix $c_{ij}$  are interpreted as a winning, provided that player $X$ has chosen the $i$-th strategy and player $Y$ has chosen the $j$-th strategy.

Let consider the step of algorithm. The prox-function is $d(x) = \sum_{i=1}^n x_i \log x_i$ (entropy) and $V_x(y) = \sum_{i=1}^n x_i \log \nicefrac{x_i}{y_i}$ (KL  divergence). The result of the proximal operator is
\begin{equation*}
    u = \text{prox}_{z_k}(\gamma_k \text{grad}(z_{k}, e_{k}, \tau, \xi_k)) = z_k \exp(-\gamma_k \text{grad}(z_{k}, e_{k}, \tau, \xi_k)),
\end{equation*}by this entry we mean: 
\begin{equation*}
    u_i = [z_k]_i \exp(-\gamma_k [\text{grad}(z_{k}, e_{k}, \tau, \xi_k)]_i).
\end{equation*} 
Using the Bregman projection onto the simplex in following way $P(x) = \nicefrac{x}{\|x\|_1}$, we have
\begin{eqnarray*}
    [x_{k+1}]_i = \frac{[x_k]_i \exp(-\gamma_k [\text{grad}_x(z_{k}, e_{k}, \tau, \xi_k)]_i)}{\sum\limits_{j=1}^n [x_k]_j \exp(-\gamma_k [\text{grad}_x(z_{k}, e_{k}, \tau, \xi_k)]_j)},
\end{eqnarray*}
\begin{eqnarray*}
    [y_{k+1}]_i = \frac{[y_k]_i \exp(\gamma_k [\text{grad}_y(z_{k}, e_{k}, \tau, \xi_k)]_i)}{\sum\limits_{j=1}^n [y_k]_j \exp(\gamma_k [\text{grad}_y(z_{k}, e_{k}, \tau, \xi_k)]_j)},
\end{eqnarray*}
where under $g_x, g_y$ we mean parts of $g$ which are responsible for $x$ and for $y$.

\begin{figure}[h!]
\centering
\begin{minipage}{0.95\textwidth}
\includegraphics[width =  \textwidth]{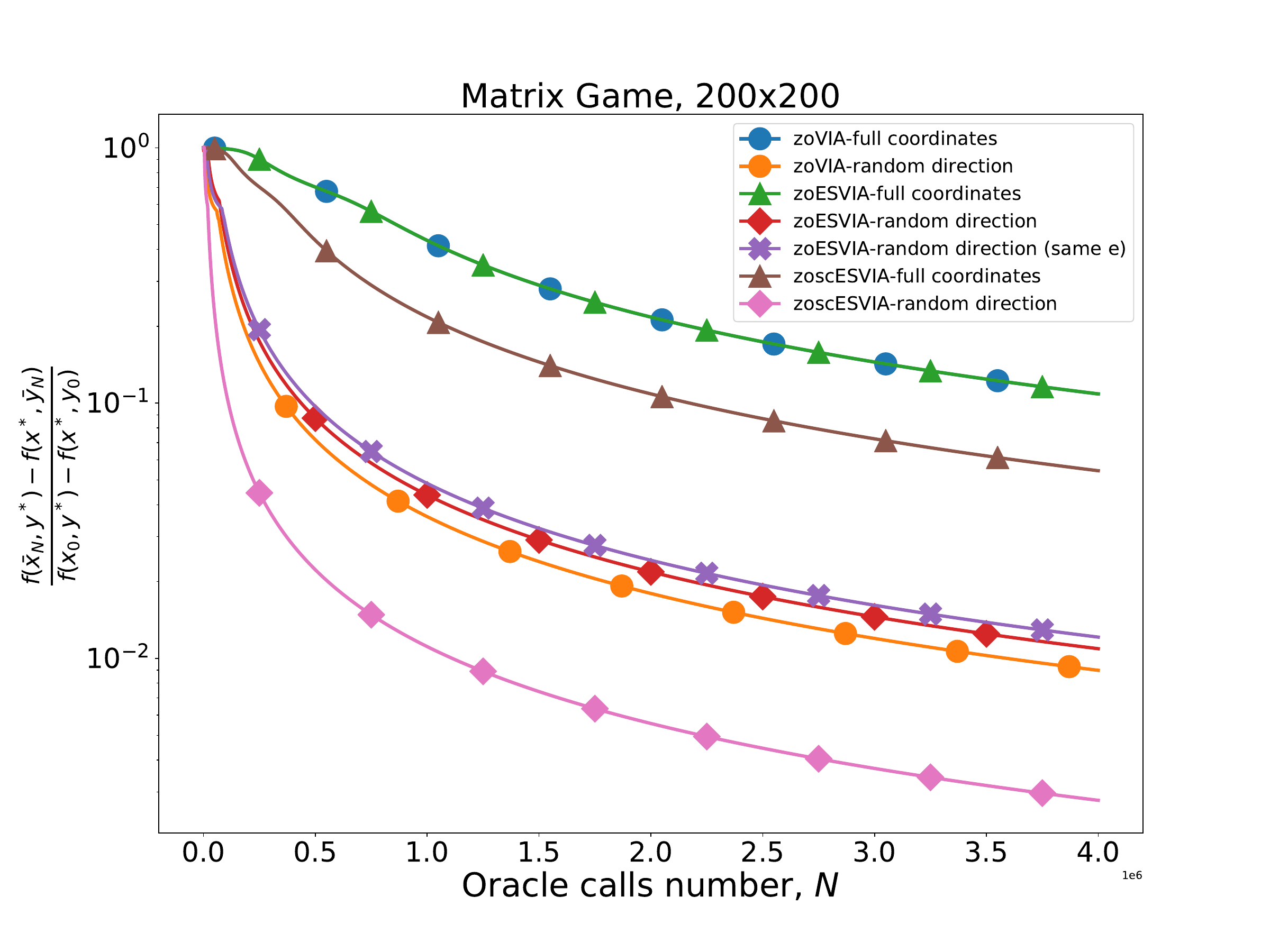}
\end{minipage}%
\caption{Different algorithms with Full coordinate and Random direction oracles applied to solve saddle-problem \eqref{exp_pr_4}. }
\label{fig:3}
\end{figure}

In the first part of the experiment, we take matrix $200 \times 200$. All elements of the matrix are generated from the uniform distribution from 0 to 1. Next, we select one row of the matrix and generate its elements from the uniform from 5 to 10. Finally, we take one element from this row and generate it uniformly from 1 to 5. The results of the experiment is on Figure \ref{fig:3}. 

From the experiment results, one can easily see the best approach in terms of oracle complexity.

\section{Conclusion}

In this paper, we presented various algorithms for optimizing smooth stochastic saddle point problems using zero-order oracles. For some oracles, we provide a theoretical analysis. We also compare the approaches covered in the work on a practical matrix game.

As a continuation of the work, we can distinguish the following areas: convergence estimates for Algorithm 4 (see the appendix), the study of gradient-free methods for saddle point problems already with a one-point approximation (in this work, we used a two-point one). We also highlight the acceleration of these methods.

\bibliographystyle{splncs04}
\bibliography{literature}

\begin{thebibliography}{10}
\providecommand{\url}[1]{\texttt{#1}}
\providecommand{\urlprefix}{URL }
\providecommand{\doi}[1]{https://doi.org/#1}

\bibitem{agarwal2010optimal}
Agarwal, A., Dekel, O., Xiao, L.: Optimal algorithms for online convex
  optimization with multi-point bandit feedback. In: COLT 2010 - The 23rd
  Conference on Learning Theory (2010)

\bibitem{basar1998dynamic}
Basar, T., Olsder, G.J.: Dynamic Noncooperative Game Theory, 2nd Edition.
  Society for Industrial and Applied Mathematics (1998).
  \doi{10.1137/1.9781611971132},
  \url{https://epubs.siam.org/doi/abs/10.1137/1.9781611971132}

\bibitem{nemirovski}
Ben-Tal, A., Nemirovski, A.: Lectures on Modern Convex Optimization: Analysis,
  Algorithms, and Engineering Applications (2019)

\bibitem{aleks2019derivativefree}
Beznosikov, A., Gorbunov, E., Gasnikov, A.: Derivative-free method for
  decentralized distributed non-smooth optimization. arXiv preprint
  arXiv:1911.10645  (2019)

\bibitem{beznosikov_sadiev_gasnikov}
Beznosikov, A., Sadiev, A., Gasnikov, A.: Gradient-free methods for
  saddle-point problem. arXiv preprint arXiv:2005.05913  (2020)

\bibitem{brent1973algorithms}
Brent, R.: Algorithms for Minimization Without Derivatives. Dover Books on
  Mathematics, Dover Publications (1973)

\bibitem{bubeck2012regret}
Bubeck, S., Cesa-Bianchi, N.: Regret analysis of stochastic and nonstochastic
  multi-armed bandit problems. Foundations and Trends® in Machine Learning
  \textbf{5}(1),  1--122 (2012). \doi{10.1561/2200000024}

\bibitem{chambolle2011first-order}
Chambolle, A., Pock, T.: A first-order primal-dual algorithm for convex
  problems with applications to imaging. Journal of Mathematical Imaging and
  Vision  \textbf{40}(1),  120--145 (2011)

\bibitem{Chen_2017}
Chen, P.Y., Zhang, H., Sharma, Y., Yi, J., Hsieh, C.J.: Zoo. Proceedings of the
  10th ACM Workshop on Artificial Intelligence and Security - AISec ’17
  (2017). \doi{10.1145/3128572.3140448},
  \url{http://dx.doi.org/10.1145/3128572.3140448}

\bibitem{choromanski2018structured}
Choromanski, K., Rowland, M., Sindhwani, V., Turner, R., Weller, A.: Structured
  evolution with compact architectures for scalable policy optimization. In:
  Dy, J., Krause, A. (eds.) Proceedings of the 35th International Conference on
  Machine Learning. Proceedings of Machine Learning Research, vol.~80, pp.
  970--978. PMLR, Stockholmsmässan, Stockholm Sweden (10--15 Jul 2018)

\bibitem{conn2009introduction}
Conn, A.R., Scheinberg, K., Vicente, L.N.: Introduction to Derivative-Free
  Optimization. Society for Industrial and Applied Mathematics (2009).
  \doi{10.1137/1.9780898718768}

\bibitem{croce2018randomized}
Croce, F., Hein, M.: A randomized gradient-free attack on relu networks (2018)

\bibitem{croce2019scaling}
Croce, F., Rauber, J., Hein, M.: Scaling up the randomized gradient-free
  adversarial attack reveals overestimation of robustness using established
  attacks (2019)

\bibitem{duchi2015optimal}
Duchi, J.C., Jordan, M.I., Wainwright, M.J., Wibisono, A.: Optimal rates for
  zero-order convex optimization: The power of two function evaluations. {IEEE}
  Trans. Information Theory  \textbf{61}(5),  2788--2806 (2015),
  arXiv:1312.2139

\bibitem{dvurechensky2020accelerated}
Dvurechensky, P., Gorbunov, E., Gasnikov, A.: An accelerated directional
  derivative method for smooth stochastic convex optimization. European Journal
  of Operational Research  (2020).
  \doi{https://doi.org/10.1016/j.ejor.2020.08.027}

\bibitem{fabian1967stochastic}
Fabian, V.: Stochastic approximation of minima with improved asymptotic speed.
  Ann. Math. Statist.  \textbf{38}(1),  191--200 (02 1967).
  \doi{10.1214/aoms/1177699070}

\bibitem{facchinei2007finite}
Facchinei, F., Pang, J.S.: Finite-dimensional variational inequalities and
  complementarity problems. Springer Science \& Business Media (2007)

\bibitem{fazel2018global}
Fazel, M., Ge, R., Kakade, S., Mesbahi, M.: Global convergence of policy
  gradient methods for the linear quadratic regulator. In: Dy, J., Krause, A.
  (eds.) Proceedings of the 35th International Conference on Machine Learning.
  Proceedings of Machine Learning Research, vol.~80, pp. 1467--1476. PMLR,
  Stockholmsmässan, Stockholm Sweden (10--15 Jul 2018)

\bibitem{fu2015handbook}
Fu, M.C., et~al.: Handbook of simulation optimization, vol.~216. Springer
  (2015)

\bibitem{gasnikov2016gradient-free}
Gasnikov, A.V., Lagunovskaya, A.A., Usmanova, I.N., Fedorenko, F.A.:
  Gradient-free proximal methods with inexact oracle for convex stochastic
  nonsmooth optimization problems on the simplex. Automation and Remote Control
   \textbf{77}(11),  2018--2034 (Nov 2016). \doi{10.1134/S0005117916110114},
  arXiv:1412.3890

\bibitem{ghadimi2013stochastic}
Ghadimi, S., Lan, G.: Stochastic first- and zeroth-order methods for nonconvex
  stochastic programming. SIAM Journal on Optimization  \textbf{23}(4),
  2341--2368 (2013), arXiv:1309.5549

\bibitem{ghadimi2016mini-batch}
Ghadimi, S., Lan, G., Zhang, H.: Mini-batch stochastic approximation methods
  for nonconvex stochastic composite optimization. Mathematical Programming
  \textbf{155}(1),  267--305 (2016). \doi{10.1007/s10107-014-0846-1},
  arXiv:1308.6594

\bibitem{gidel2018variational}
Gidel, G., Berard, H., Vignoud, G., Vincent, P., Lacoste-Julien, S.: A
  variational inequality perspective on generative adversarial networks (2018)

\bibitem{goodfellow2014generative}
Goodfellow, I.J., Pouget-Abadie, J., Mirza, M., Xu, B., Warde-Farley, D.,
  Ozair, S., Courville, A., Bengio, Y.: Generative adversarial networks (2014)

\bibitem{goodfellow2014explaining}
Goodfellow, I.J., Shlens, J., Szegedy, C.: Explaining and harnessing
  adversarial examples (2014)

\bibitem{gorbunov2018accelerated}
Gorbunov, E., Dvurechensky, P., Gasnikov, A.: An accelerated method for
  derivative-free smooth stochastic convex optimization. arXiv preprint
  arXiv:1802.09022  (2018)

\bibitem{harker1990finite}
Harker, P.T., Pang, J.S.: Finite-dimensional variational inequality and
  nonlinear complementarity problems: a survey of theory, algorithms and
  applications. Mathematical programming  \textbf{48}(1-3),  161--220 (1990)

\bibitem{hsieh2019convergence}
Hsieh, Y.G., Iutzeler, F., Malick, J., Mertikopoulos, P.: On the convergence of
  single-call stochastic extra-gradient methods (2019)

\bibitem{juditsky2008solving}
Juditsky, A., Nemirovskii, A.S., Tauvel, C.: Solving variational inequalities
  with stochastic mirror-prox algorithm (2008)

\bibitem{Korpelevich1976TheEM}
Korpelevich, G.M.: The extragradient method for finding saddle points and other
  problems (1976)

\bibitem{liu2019minmax}
Liu, S., Lu, S., Chen, X., Feng, Y., Xu, K., Al-Dujaili, A., Hong, M.,
  O'Reilly, U.M.: Min-max optimization without gradients: Convergence and
  applications to adversarial ml (2019)

\bibitem{madry2018towards}
Madry, A., Makelov, A., Schmidt, L., Tsipras, D., Vladu, A.: Towards deep
  learning models resistant to adversarial attacks. In: 6th International
  Conference on Learning Representations, {ICLR} 2018, Vancouver, BC, Canada,
  April 30 - May 3, 2018, Conference Track Proceedings (2018)

\bibitem{morgenstern1953theory}
Morgenstern, O., Von~Neumann, J.: Theory of games and economic behavior.
  Princeton university press (1953)

\bibitem{conf/cvpr/NarodytskaK17}
Narodytska, N., Kasiviswanathan, S.P.: Simple black-box adversarial attacks on
  deep neural networks. In: CVPR Workshops. pp. 1310--1318. IEEE Computer
  Society (2017),
  \url{http://doi.ieeecomputersociety.org/10.1109/CVPRW.2017.172}

\bibitem{nedic2009subgradient}
Nedi{\'{c}}, A., Ozdaglar, A.: Subgradient methods for saddle-point problems.
  Journal of Optimization Theory and Applications  \textbf{142}(1),  205--228
  (Jul 2009)

\bibitem{Nemirovski2004}
Nemirovski, A.: Prox-method with rate of convergence o (1/ t ) for variational
  inequalities with lipschitz continuous monotone operators and smooth
  convex-concave saddle point problems. SIAM Journal on Optimization
  \textbf{15},  229--251 (01 2004). \doi{10.1137/S1052623403425629}

\bibitem{nesterov2017random}
Nesterov, Y., Spokoiny, V.: Random gradient-free minimization of convex
  functions. Found. Comput. Math.  \textbf{17}(2),  527--566 (Apr 2017).
  \doi{10.1007/s10208-015-9296-2}, first appeared in 2011 as CORE discussion
  paper 2011/16

\bibitem{pinto2017robust}
Pinto, L., Davidson, J., Sukthankar, R., Gupta, A.: Robust adversarial
  reinforcement learning. Proceedings of Machine Learning Research, vol.~70,
  pp. 2817--2826. PMLR, International Convention Centre, Sydney, Australia
  (06--11 Aug 2017), \url{http://proceedings.mlr.press/v70/pinto17a.html}

\bibitem{rosenbrock1960automatic}
Rosenbrock, H.H.: An automatic method for finding the greatest or least value
  of a function. The Computer Journal  \textbf{3}(3),  175--184 (1960).
  \doi{10.1093/comjnl/3.3.175}

\bibitem{salimans2017evolution}
Salimans, T., Ho, J., Chen, X., Sidor, S., Sutskever, I.: Evolution strategies
  as a scalable alternative to reinforcement learning. arXiv:1703.03864  (2017)

\bibitem{shamir2017optimal}
Shamir, O.: An optimal algorithm for bandit and zero-order convex optimization
  with two-point feedback. Journal of Machine Learning Research  \textbf{18},
  52:1--52:11 (2017)

\bibitem{shashaani2018ASTRO}
Shashaani, S., Hashemi, F.S., Pasupathy, R.: Astro-df: A class of adaptive
  sampling trust-region algorithms for derivative-free stochastic optimization.
  SIAM Journal on Optimization  \textbf{28}(4),  3145--3176 (2018).
  \doi{10.1137/15M1042425}, \url{https://doi.org/10.1137/15M1042425}

\bibitem{spall2003introduction}
Spall, J.C.: Introduction to Stochastic Search and Optimization. John Wiley \&
  Sons, Inc., New York, NY, USA, 1 edn. (2003)

\bibitem{stich2019unified}
Stich, S.U.: Unified optimal analysis of the (stochastic) gradient method
  (2019)

\bibitem{stich2011linear}
Stich, S.U., Muller, C.L., Gartner, B.: Optimization of convex functions with
  random pursuit. SIAM Journal on Optimization  \textbf{23}(2),  1284--1309
  (2013)

\bibitem{tramr2017ensemble}
Tramèr, F., Kurakin, A., Papernot, N., Goodfellow, I., Boneh, D., McDaniel,
  P.: Ensemble adversarial training: Attacks and defenses (2017)

\bibitem{wang2020zerothorder}
Wang, Z., Balasubramanian, K., Ma, S., Razaviyayn, M.: Zeroth-order algorithms
  for nonconvex minimax problems with improved complexities (2020)

\bibitem{ye2018hessianaware}
Ye, H., Huang, Z., Fang, C., Li, C.J., Zhang, T.: Hessian-aware zeroth-order
  optimization for black-box adversarial attack (2018)

\end{thebibliography}

\clearpage

\appendix

\section{General facts and technical lemmas}\label{sec:gf_tl}

\begin{lemma}
For arbitrary integer $n\ge 1$ and arbitrary set of positive numbers $a_1,\ldots,a_n$ we have
\begin{equation}
    \left(\sum\limits_{i=1}^m a_i\right)^2 \le m\sum\limits_{i=1}^m a_i^2.\label{eq:squared_sum}
\end{equation}
\end{lemma}

\begin{lemma}
For $q\ge 2$ and for arbitrary vectors $a \in \RR^n, b \in \RR^m$ we have
\begin{equation}
    \left\|\begin{pmatrix}
a\\
b
\end{pmatrix} \right\|^2_q \le \| a\|^2_q + \| b\|^2_q.\label{eq:block_vec}
\end{equation}
\end{lemma}

\begin{lemma}[Fact 5.3.2 from \cite{nemirovski}]
    Given norm $\|\cdot \|$ on space $\mathcal{Z}$ and prox-function $d(z)$, let $z \in \mathcal{Z}$, $w \in \mathbb{R}^n$ and $z_{+} = \text{prox}_z(w)$. Then for all $u \in \mathcal{Z}$
    \begin{eqnarray}
        \label{lemma3_1}
         \langle w, z_{+} - u\rangle \leqslant V_{z}(u) - V_{z_{+}}(u) - V_{z}(z_{+}).
    \end{eqnarray}
\end{lemma}

\begin{lemma} [see Lemma 1 from \cite{gorbunov2018accelerated}]
Let $e \in \mathcal{RS}^2(1)$, i.e. a random vector uniformly distributed on the surface of the unit Euclidean sphere in $\RR^n$, $q \in [2;+\infty)$. Then, for $n \geq 8$, 
\begin{eqnarray}
\label{tlemma_1}
\EE\left[\|e\|^2_q\right] &\leq&  n^{2/q - 1} \rho_n, \\
\label{tlemma_2}
\EE\left[ \langle s,e \rangle^2 \|e \|^2_q\right] &\leq& 6 n^{2/q - 2} \rho_n \|s\|^2_2, ~~\forall s \in \RR^n,
\end{eqnarray}
where $\rho_n = \min\{q-1, 16 \log n -8\}$.
\end{lemma}

\begin{lemma}[see Lemma 2 from \cite{stich2019unified}] Let consider non-negative sequence $r_k$:
$$r_{k+1} \leq (1 - a \gamma)r_k + c \gamma^2,$$
where $a,c > 0$, $\gamma = \min\left(\frac{1}{d}, \frac{\log(\max(2, a^2 r_0 T / c)}{aN}\right)$. Then 
\begin{equation}
    \label{stich}
   r_{N+1} \leq r_0 \cdot \exp\left(- \frac{aN}{2d}\right) + \frac{c}{a^2 N}.
\end{equation}
\end{lemma}

\section{Proof of Lemma \ref{lemma2}} 

\textbf{Lemma.} If $f(x,y)$ is $\mu$-strongly convex on $x$ and $\mu$-strongly concave on $y$ w.r.t $V_{\cdot}(\cdot)$, then for $F(z)$ we have
\begin{equation}
    \label{lemma2_a}
    \langle F(z_1) - F(z_2), z_1 - z_2 \rangle \geq \frac{\mu}{2}\left( V_{z_1}(z_2) + V_{z_2}(z_1)\right),~~ \forall z_1, z_2 \in \mathcal{Z}.
\end{equation}

\textbf{Proof.}
By definition of $\mu$-strong convexity w.r.t $V_{\cdot}(\cdot)$:
$$f(x_1,y_2) \geq f(x_2, y_2) + \langle\nabla_x f(x_2, y_2) , x_1 - x_2\rangle + \frac{\mu}{2}\left(V_{(x_2, y_2)}(x_1, y_2) + V_{(x_1, y_2)}(x_2, y_2)\right) ,$$
$$f(x_2,y_1) \geq f(x_1, y_1) + \langle\nabla_x f(x_1, y_1) , x_2 - x_1\rangle + \frac{\mu}{2}\left(V_{(x_1, y_1)}(x_2, y_1) + V_{(x_2, y_1)}(x_1, y_1)\right),$$
$$-f(x_1,y_2) \geq -f(x_1, y_1) + \langle-\nabla_y f(x_1, y_1) , y_2 - y_1\rangle + \frac{\mu}{2}\left(V_{(x_1, y_1)}(x_1, y_2) + V_{(x_1, y_2)}(x_1, y_1)\right),$$
$$-f(x_2,y_1) \geq -f(x_2, y_2) + \langle-\nabla_y f(x_2, y_2) , y_1 - y_2\rangle + \frac{\mu}{2}\left(V_{(x_2, y_2)}(x_2, y_1) + V_{(x_1, y_1)}(x_2, y_2)\right) .$$
Let introduce a new definition for sum of Bregman divergences:
\begin{eqnarray*} 
\mathcal{V} &=& V_{(x_2, y_2)}(x_1, y_2) + V_{(x_1, y_2)}(x_2, y_2) + V_{(x_1, y_1)}(x_2, y_1) + V_{(x_2, y_1)}(x_1, y_1)\\&& + V_{(x_1, y_1)}(x_1, y_2) + V_{(x_1, y_2)}(x_1, y_1) + V_{(x_2, y_2}(x_2, y_1) + V_{(x_1, y_1)}(x_2, y_2).
\end{eqnarray*}
Using definition of Bregman divergence and $1$-stronge convexity of prox-function $d$, we get:
\begin{eqnarray*} 
\mathcal{V} &=& \langle\nabla_x d(x_2, y_2) - \nabla_x d(x_1, y_2), x_2 - x_1\rangle \\&&+ \langle\nabla_x d(x_2, y_1) - \nabla_x d(x_1, y_1), x_2 - x_1\rangle\\&& + \langle\nabla_y d(x_2, y_2) - \nabla_y d(x_2, y_1), y_2 - y_1\rangle\\&& + \langle\nabla_y d(x_1, y_2) - \nabla_y d(x_1, y_1), y_2 - y_1\rangle \\ &=& \langle\nabla d(z_2) - \nabla d(z_1), z_2 - z_1\rangle  + \langle\nabla d(\tilde{z}_2) - \nabla d(\tilde{z}_1), \tilde{z}_2 - \tilde{z}_1\rangle \\ &\geq& V_{z_1}(z_2) + V_{z_2}(z_1), 
\end{eqnarray*} where $\tilde{z}_2 =(x_2, y_1) $, $\tilde{z}_1 =(x_1, y_2) $
Thus, we have $\mathcal{V} \geq V_{z_1}(z_2) + V_{z_2}(z_1)$.
Summming up:
\begin{eqnarray*} 
\langle\nabla_x f(x_2, y_2) - \nabla_x f(x_1, y_1) , x_1 - x_2\rangle \quad \quad \quad \quad \quad \quad \quad \quad \quad \quad \quad \quad \\
- \langle\nabla_y f(x_2, y_2) - \nabla_y f(x_1, y_1) , y_1 - y_2\rangle + \frac{\mu \mathcal{V}}{2}\leq 0.
\end{eqnarray*}
Using $\mathcal{V} \geq V_{z_1}(z_2) + V_{z_2}(z_1)$, we have
\begin{eqnarray*} 
\langle\nabla_x f(x_2, y_2) - \nabla_x f(x_1, y_1) , x_1 - x_2\rangle - \langle\nabla_y f(x_2, y_2) - \nabla_y f(x_1, y_1) , y_1 - y_2\rangle\\ + \frac{\mu}{2}(V_{z_1}(z_2) + V_{z_2}(z_1))\leq 0, 
\end{eqnarray*}
and 
\begin{eqnarray*}
\langle F(z_1) - F(z_2), z_1 - z_2 \rangle &=& \langle\nabla_x f(x_2, y_2) - \nabla_x f(x_1, y_1) , x_2 - x_1\rangle\\ &&-\langle\nabla_y f(x_2, y_2) - \nabla_y f(x_1, y_1) , y_2 - y_1\rangle\\ &\geq& \frac{\mu}{2}(V_{z_1}(z_2) + V_{z_2}(z_1)) .
\end{eqnarray*}
\EndProof

\section{Proof of Lemma \ref{lemma1}}

\textbf{Lemma.}
Let $e \in \mathcal{RS}^2(1)$, i.e. uniformly distributed on the unit Euclidean sphere. Randomness comes from independent variables $e$, $\xi$ and a point $z$. Norm $\|\cdot\|_* = \|\cdot\|_q$ satisfies $q \in [2;+\infty)$. We introduce the constant $\rho_n$:
\begin{equation*}
    \rho_n = \min\{q-1, 16 \log(n) -8 \}.
\end{equation*}
Then under Assumption 3 or 3(f) the following statements hold:
\begin{itemize}
    \item for Random direction oracle
    \begin{eqnarray}
        \label{lemma1_2_a}
        \EE\left[\|g_d(z,e,\tau, \xi)\|^2_q\right] &\leq& 48 n^{2/q} \rho_n \EE\left[\| F(z) -F(z^*) \|^2_2\right] + 48 n^{2/q} \rho_n \| F(z^*) \|^2_2 \nonumber\\ 
    && + 48 n^{2/q} \rho_n \sigma^2 + 8 n^{2/q + 1}\rho_n L^2\tau^2  \nonumber \\
    && + 16 \frac{n^{2/q + 1}\rho_n \Delta^2}{\tau^2}, \\
        \label{lemma1_2_2_a}
        \left\| \EE [g_d(z,e,\tau, \xi)] - F(z)\right\|_q &\leq& 2 n^{1/q + 1/2} \sqrt{\rho_n} L\tau + 4 n^{1/q + 1/2} \sqrt{\rho_n} \frac{\Delta}{\tau};
    \end{eqnarray}
    \item for Full coordinates oracle 
    \begin{eqnarray}
        \label{lemma1_3_a}
        \EE\left[\|g_f(z,\tau, \xi) - F(z)\|^2_q\right] &\leq& 3\sigma^2 + 3nL_2^2 \tau^2 + \frac{6n\Delta^2}{\tau^2}, \\
        \label{lemma1_3_1_a}
        \|\EE\left[g_f(z,\tau, \xi)] - F(z)\right\|_q &\leq& \sqrt{n} L \tau + \frac{2\sqrt{n}\Delta}{\tau}.
    \end{eqnarray}
\end{itemize}

\textbf{Proof of} \eqref{lemma1_2_a}.

\begin{eqnarray*} 
\EE\left[\|g_d(z,e,\tau, \xi)\|^2_q\right] &\overset{\eqref{eq:squared_sum}}{\leq}& 4 n^2 \EE\left[\left\|\left(\begin{array}{c}
    \left\langle \nabla_x f(x,y) , e_x \right\rangle e_x\\
    \left\langle -\nabla_y f(x,y) , e_y \right\rangle e_y 
    \end{array}\right)\right\|^2_q\right] \nonumber\\
&&+ 4 n^2\EE\left[\left\| \left(\begin{array}{c}
    \left\langle \nabla_x f(x,y, \xi) - \nabla_x f(x,y) , e_x \right\rangle e_x\\
    \left\langle -\nabla_y f(x,y, \xi) +\nabla_y f(x,y) , e_y \right\rangle e_y 
    \end{array}\right)\right\|^2_q\right] \nonumber\\
&&+ 4 \frac{n^2}{\tau^2} \EE\left[ \left\| \left(
    \begin{array}{c}
    \left( f(x + \tau e_x, y,  \xi) -   f(x, y, \xi) - \langle \nabla_x f(x,y, \xi), \tau e_x\rangle  \right)e_x\\
    \left( f(x, y, \xi) - f(x, y + \tau e_y,  \xi) +  \langle \nabla_y f(x,y, \xi), \tau e_y\rangle  \right)e_y 
    \end{array}
    \right) \right\|^2_q\right] \nonumber
\\&&+ 4 \frac{n^2}{\tau^2} \EE\left[ \left\| \left(
    \begin{array}{c}
    \left(\delta(x + \tau e_x, y) -  \delta(x, y) \right)e_x\\
    \left(\delta(x, y) - \delta(x, y + \tau e_y) \right)e_y 
    \end{array}
    \right) \right\|^2_q\right] \nonumber\\
&\overset{\eqref{eq:block_vec}}{\leq}&   4 n^2 \EE\left[\left\|
    \left\langle \nabla_x f(x,y) ,  e_x \right\rangle e_x\right\|^2_q\right] +
    4 n^2 \EE\left[\left\|\left\langle -\nabla_y f(x,y) , e_y \right\rangle e_y 
    \right\|^2_q\right] \nonumber\\
&&+ 4 n^2\EE\left[\left\| 
    \left\langle \nabla_x f(x,y, \xi) - \nabla_x f(x,y) ,  e_x \right\rangle e_x\right\|^2_q\right] \\
    &&+ 4 n^2\EE\left[\left\| \left\langle -\nabla_y f(x,y, \xi) +\nabla_y f(x,y) ,  e_y \right\rangle e_y 
    \right\|^2_q\right] \nonumber\\
&&+ 4 \frac{n^2}{\tau^2} \EE\left[ \left\| 
    \left(\tilde f(x + \tau e_x, y,  \xi) -  \tilde f(x, y, \xi) - \langle \nabla_x f(x,y, \xi), \tau e_x\rangle  \right)e_x  \right\|^2_q\right]\\
    &&+ 4 \frac{n^2}{\tau^2} \EE\left[ \left\| \left(\tilde f(x, y, \xi) - \tilde f(x, y + \tau e_y,  \xi) +  \langle \nabla_y f(x,y, \xi), \tau e_y\rangle  \right)e_y 
     \right\|^2_q\right] \nonumber
\\&&+ 4 \frac{n^2}{\tau^2} \EE\left[ \left\| 
    \left(\delta(x + \tau e_x, y) -  \delta(x, y) \right)e_x\right\|^2_q\right]\\
  &&+ 4 \frac{n^2}{\tau^2} \EE\left[ \left\| \left(\delta(x, y) - \delta(x, y + \tau e_y) \right)e_y 
     \right\|^2_q\right].
\end{eqnarray*}
From \eqref{SP_c} we get $\|\nabla_x f(x_1, y, \xi) - \nabla_x f(x_2, y, \xi) \|_2 \leq L \|x_1 - x_2 \|_2$ and $\|\nabla_y f(x, y_1, \xi) - \nabla_y f(x, y_2, \xi) \|_2 \leq L \|y_1 - y_2 \|_2$ for all $x, x_1, x_2 \in \mathcal{X}$, $y, y_1, y_2 \in \mathcal{Y}$. It follows that functions $f(\cdot, y, \xi)$ and $f(x, \cdot, \xi)$ are $L(\xi)$-Lipschitz continuous. Then 

\begin{eqnarray*} 
\EE\left[\|g_d(z,e,\tau, \xi)\|^2_q\right] &\overset{}{\leq}&   4 n^2 \EE\left[\left\|
    \left\langle \nabla_x f(x,y) , \tau e_x \right\rangle e_x\right\|^2_q\right] +
    4 n^2 \EE\left[\left\|\left\langle -\nabla_y f(x,y) , \tau e_y \right\rangle e_y 
    \right\|^2_q\right] \nonumber\\
&&+ 4 n^2\EE\left[\left\| 
    \left\langle \nabla_x f(x,y, \xi) - \nabla_x f(x,y) , \tau e_x \right\rangle e_x\right\|^2_q\right] \\
    &&+ 4 n^2\EE\left[\left\| \left\langle -\nabla_y f(x,y, \xi) +\nabla_y f(x,y) , \tau e_y \right\rangle e_y 
    \right\|^2_q\right] \nonumber\\
&&+ 4 n^2 L_2^2\tau^2  
    \EE\left[ \left\| e_x  \right\|^2_q\right] + 4 n^2 L_2^2\tau^2  
    \EE\left[ \left\| e_y  \right\|^2_q\right] \nonumber
\\&&+ 8 \frac{n^2 \Delta^2}{\tau^2} \EE\left[ \left\| e_x\right\|^2_q\right]+ 8 \frac{n^2 \Delta^2}{\tau^2} \EE\left[ \left\| e_y\right\|^2_q\right].
\end{eqnarray*}
In the last inequality, we additionally use \eqref{ran_dir1} + \eqref{eq:squared_sum} and 
independence of $e$ and $\xi$. With \eqref{tlemma_1} and \eqref{tlemma_2}, one can get the following result:
\begin{eqnarray*} 
\EE\left[\|g_d(z,e,\tau, \xi)\|^2_q\right] &\overset{}{\leq}&   24 n^{2/q} \rho_n \EE\left[\|
     \nabla_x f(x,y) \|^2_2\right] + 24 n^{2/q} \rho_n \EE\left[\| -\nabla_y f(x,y) \|^2_2\right] \nonumber\\
&&+ 24 n^{2/q} \rho_n\EE\left[\left\| 
    \nabla_x f(x,y, \xi) - \nabla_x f(x,y) \right\|^2_2\right] \\
    &&+ 24 n^{2/q} \rho_n\EE\left[\left\|  -\nabla_y f(x,y, \xi) +\nabla_y f(x,y)\right\|^2_2\right] \nonumber\\
&&+ 8 n^{2/q + 1}\rho_n L^2\tau^2  
     + 16 \frac{n^{2/q + 1}\rho_n \Delta^2}{\tau^2} \\
     &\overset{\eqref{ran_dir1}}{\leq}& 24 n^{2/q} \rho_n \EE\left[\| F(z) \|^2_2\right] + 48 n^{2/q} \rho_n \sigma^2 + 8 n^{2/q + 1}\rho_n L_2^2\tau^2  
     + 16 \frac{n^{2/q + 1}\rho_n \Delta^2}{\tau^2} \\
     &\overset{\eqref{eq:squared_sum}}{\leq}& 48 n^{2/q} \rho_n \EE\left[\| F(z) -F(z^*) \|^2_2\right] + 48 n^{2/q} \rho_n \| F(z^*) \|^2_2 \\ 
    && + 48 n^{2/q} \rho_n \sigma^2 + 8 n^{2/q + 1}\rho_n L_2^2\tau^2  
     + 16 \frac{n^{2/q + 1}\rho_n \Delta^2}{\tau^2}.
\end{eqnarray*}

\textbf{Proof of} \eqref{lemma1_2_2_a} .
\begin{eqnarray*} 
\left\| \EE [g_d(z,e,\tau, \xi)] - F(z)\right\|_q &\leq& \frac{n}{\tau} \left\|\EE\left[ \left(
    \begin{array}{c}
    \left( f(x + \tau e_x, y,  \xi) -   f(x, y, \xi) - \langle \nabla_x f(x,y, \xi), \tau e_x\rangle  \right)e_x\\
    \left( f(x, y, \xi) - f(x, y + \tau e_y,  \xi) +  \langle \nabla_y f(x,y, \xi), \tau e_y\rangle  \right)e_y 
    \end{array}
    \right) \right] \right\|_q  \\
&&+n \left\|\EE\left[\left(\begin{array}{c}
    \left\langle \nabla_x f(x,y, \xi) - \nabla_x f(x,y) , e_x \right\rangle e_x\\
    \left\langle -\nabla_y f(x,y, \xi) +\nabla_y f(x,y) , e_y \right\rangle e_y 
    \end{array}\right) \right] \right\|_q  \\
&&+ n \left\|\EE\left[\left(\begin{array}{c}
    \left\langle \nabla_x f(x,y) , e_x \right\rangle e_x\\
    \left\langle -\nabla_y f(x,y) , e_y \right\rangle e_y 
    \end{array}\right) \right] - \left(\begin{array}{c}
    \nabla_x f(x,y) \\
    -\nabla_y f(x,y)
    \end{array}\right)\right\|_q \nonumber\\
&&+ \frac{n}{\tau} \left\|\EE\left[ \left(
    \begin{array}{c}
    \left(\delta(x + \tau e_x, y) -  \delta(x, y) \right)e_x\\
    \left(\delta(x, y) - \delta(x, y + \tau e_y) \right)e_y 
    \end{array}
    \right)\right]\right\|_q.
\end{eqnarray*}
Taking into account the independence of $e$ and $\xi$, as well as using their unbiasedness, we get
\begin{eqnarray*} 
\left\| \EE [g_d(z,e,\tau, \xi)] - F(z)\right\|_q &\leq& \frac{n}{\tau} \left\|\EE\left[ \left(
    \begin{array}{c}
    \left( f(x + \tau e_x, y) -   f(x, y) - \langle \nabla_x f(x,y), \tau e_x\rangle  \right)e_x\\
    \left( f(x, y) - f(x, y + \tau e_y) +  \langle \nabla_y f(x,y), \tau e_y\rangle  \right)e_y 
    \end{array}
    \right) \right] \right\|_q   \\
&&+ \frac{n}{\tau} \left\|\EE\left[ \left(
    \begin{array}{c}
    \left(\delta(x + \tau e_x, y) -  \delta(x, y) \right)e_x\\
    \left(\delta(x, y) - \delta(x, y + \tau e_y) \right)e_y 
    \end{array}
    \right)\right]\right\|_q \\
    &\overset{\eqref{eq:block_vec}}{\leq}&  
    \frac{n}{\tau} \left\|\EE\left[ \left( f(x + \tau e_x, y) -   f(x, y) - \langle \nabla_x f(x,y), \tau e_x\rangle  \right)e_x  \right] \right\|_q \\ &&+ \frac{n}{\tau} \left\|\EE\left[ \left( f(x, y) - f(x, y + \tau e_y) +  \langle \nabla_y f(x,y), \tau e_y\rangle  \right)e_y 
     \right] \right\|_q   \\
&&+ \frac{n}{\tau} \left\|\EE\left[ \left(\delta(x + \tau e_x, y) -  \delta(x, y) \right)e_x\right]\right\|_q \\ &&+
     \frac{n}{\tau} \left\|\EE\left[ \left(\delta(x, y) - \delta(x, y + \tau e_y) \right)e_y \right]\right\|_q.
\end{eqnarray*}
Further, Jensen inequality gives
\begin{eqnarray*} 
\left\| \EE [g_d(z,e,\tau, \xi)] - F(z)\right\|_q &\leq&   
    \frac{n}{\tau} \EE\left[\left| f(x + \tau e_x, y) -   f(x, y) - \langle \nabla_x f(x,y), \tau e_x\rangle  \right| \left\| e_x  \right\|_q\right]  \\ &&+ \frac{n}{\tau} \EE\left[\left| f(x, y) - f(x, y + \tau e_y) +  \langle \nabla_y f(x,y), \tau e_y\rangle  \right| \left\| e_y 
      \right\|_q\right]   \\
&&+ \frac{n}{\tau} \EE\left[ \left|\delta(x + \tau e_x, y) -  \delta(x, y) \right|\left\|e_x\right\|_q\right] \\ &&+
     \frac{n}{\tau} \EE\left[\left|\delta(x, y) - \delta(x, y + \tau e_y) \right| \left\|e_y\right\|_q \right].
\end{eqnarray*}
It remains to use $L$-Lipschitz continuous of $f(\cdot, y)$ and $f(x, \cdot)$:
\begin{eqnarray*} 
\left\| \EE [g_d(z,e,\tau, \xi)] - F(z)\right\|_q &\leq&   
    nL\tau\EE\left[ \left\| e_x  \right\|_q\right]  + nL\tau \EE\left[ \left\| e_y 
      \right\|_q\right]   \\
&&+ \frac{n}{\tau} \EE\left[ \left(|\delta(x + \tau e_x, y)| +  |\delta(x, y)| \right)\left\|e_x\right\|_q\right] \\ &&+
     \frac{n}{\tau} \EE\left[\left(|\delta(x, y)| + |\delta(x, y + \tau e_y)| \right) \left\|e_y\right\|_q \right] \\
     &\overset{\eqref{ran_dir1}, \eqref{tlemma_1}}{\leq}& 2 n^{1/q + 1/2} \sqrt{\rho_n} L\tau + 4 n^{1/q + 1/2} \sqrt{\rho_n} \frac{\Delta}{\tau}.
\end{eqnarray*}



\textbf{Proof of} \eqref{lemma1_3_a}.
\begin{eqnarray*} 
\EE\left[\|g_f(z,\tau, \xi) - F(z)\|^2_q\right] &\overset{\eqref{full_coor},\eqref{eq:squared_sum}}{\leq}& 3\EE\Bigg[\Bigg\| \frac{1}{\tau}\sum\limits_{i=1}^{n_x}\left( f(z+\tau h_i, \xi) -  f(z, \xi) \right)h_i \\
&&+ \frac{1}{\tau}\sum\limits_{i=n_x + 1}^{n_x + n_y}\left( f(z, \xi)  -  f(z+\tau h_i, \xi) \right)h_i - F(z, \xi)  \Bigg\|^2_2\Bigg] \nonumber\\
&&+ 3\EE\left[\left\| F(z, \xi)  - F(z)\right\|^2_2\right] \\
&&+ 3\EE\Bigg[\Bigg\| \sum\limits_{i=1}^{n_x+n_y}\frac{\left( \delta(z+\tau h_i)  -\delta(z)\right)}{{\tau}}h_i \Bigg\|^2_2\Bigg] \\
&\overset{\eqref{ran_dir1},\eqref{eq:squared_sum}}{\leq}&
3\EE\left[\sum\limits_{i=1}^{n_x+n_y}\left| \frac{\left( f(z+\tau h_i, \xi) -  f(z, \xi) \right)}{{\tau}} - \frac{\partial f(z, \xi)}{\partial z_i}  \right|^2\right] \nonumber\\
&&+ 3 \sigma^2 + 6 \frac{n\Delta^2}{\tau^2}.
\end{eqnarray*}

By the mean value theorem we have that for some $|q_i| \leq |\tau|$:
\begin{eqnarray*} 
\EE\left[\|g_f(z,\tau, \xi) - F(z)\|^2_*\right] 
&\overset{}{\leq}&
3\EE\left[\sum\limits_{i=1}^n\left|  \frac{\partial f(z + q_i h_i, \xi)}{\partial z_i}  - \frac{\partial f(z, \xi)}{\partial z_i}  \right|^2\right] \nonumber\\
&&+ 3 \sigma^2 + 6 \frac{n\Delta^2}{\tau^2} \\
&\overset{}{\leq}& 3\sum\limits_{i=1}^n L_2^2 q_i^2+ 3 \sigma^2 + 6 \frac{n\Delta^2}{\tau^2}
\\
&\overset{}{\leq}& 3n L_2^2 \tau^2 + 3 \sigma^2 + 6 \frac{n\Delta^2}{\tau^2}
.
\end{eqnarray*}

\textbf{Proof of} \eqref{lemma1_3_1_a}.
Using unbiasedness of $\xi$:
\begin{eqnarray*} 
\left\|\EE\left[g_f(z,\tau, \xi)] - F(z)\right] \right\|_q &\leq&  
\Bigg\| \frac{1}{\tau}\sum\limits_{i=1}^{n_x}\left( f(z+\tau h_i) -  f(z) \right)h_i \\
&&+ \frac{1}{\tau}\sum\limits_{i=n_x + 1}^{n_x + n_y}\left( f(z)  -  f(z+\tau h_i, ) \right)h_i - F(z)  \Bigg\|_2 \nonumber\\
&&+ \Bigg\| \sum\limits_{i=1}^{n_x+n_y}\frac{\left( \delta(z+\tau h_i)  -\delta(z)\right)}{{\tau}}h_i \Bigg\|_2 \\
&\overset{\eqref{ran_dir1}}{\leq}& \sqrt{\sum\limits_{i=1}^n L^2 q_i^2} + \frac{2\sqrt{n}\Delta}{\tau}
\\
&\overset{}{\leq}& \sqrt{n} L \tau + \frac{2\sqrt{n}\Delta}{\tau}.
\end{eqnarray*}
\EndProof

\section{Proof of Theorem 1} 

\begin{lemma}\label{lemmt_1}
Let $z, g\in \RR^n$ and $\mathcal{Z} \subset \RR^n$. Then for $z_{1} = \text{prox}_z(g)$ and for all $u \in \mathcal{Z}$ we have
\begin{equation}
\label{DL_a}
   \langle g , z - u\rangle \leq
    V_{z}(u) - V_{z_{1}}(u) + \frac{1}{2}\| g \|^2_q
\end{equation}
\end{lemma}

\textbf{Proof.} By \eqref{lemma3_1}, we have for all $u \in \mathcal{Z}$
\begin{eqnarray*}
    \langle g , z_{1} - u\rangle = \langle g , z_{1} - z + z - u\rangle \leq V_{z}(u) - V_{z_{1}}(u) - V_{z}(z_{1}).
\end{eqnarray*}
Making simple transformations:
\begin{eqnarray*}
    \langle g , z - u\rangle &\leq& \langle g , z - z_{1}\rangle + V_{z}(u) - V_{z_{1}}(u) - V_{z}(z_{1}) \nonumber\\
    &\leq& \langle g , z - z_{1}\rangle + V_{z}(u) - V_{z_{1}}(u) - \frac{1}{2}\|z_{1} - z\|^2_p. 
\end{eqnarray*}
In last inequality we use the property of the Bregman divergence: $V_x(y) \ge \frac{1}{2}\|x-y\|_p^2$. With Hölder's inequality and the fact: $ab - \nicefrac{b^2}{2} \leqslant\nicefrac{a^2}{2}$, we get
\begin{eqnarray*}
    \langle g, z - u\rangle &\leq&
    \| g\|_q\|z - z_{1}\|_p  + V_{z}(u) - V_{z_{1}}(u) - \frac{1}{2}\|z_{1} - z\|^2_p  \nonumber\\
    &\leq&
    V_{z}(u) - V_{z_{1}}(u) + \frac{1}{2}\| g \|^2_q.
\end{eqnarray*}
\EndProof

\textbf{Theorem.} By Algorithm 1 with Random direction oracle
\begin{itemize}
    \item under Assumptions 1, 2, 3(f) and with $\gamma \leq \frac{1}{48 n^{\nicefrac{2}{q}}\rho_nL}$, we get
    \begin{eqnarray}
    \label{teor1.1_a}
        \frac{1}{N+1}\sum^N_{k = 0}\EE\left[\|F(z_k) - F(z^*)\|^2_2\right] &\leq&\frac{2LD^2_p}{\gamma N} + 48 \gamma n^{2/q} \rho_nL\left(\|F(z^*)\|^2_2 + \sigma^2\right) \nonumber\\&& + 8\gamma n^{2/q + 1} \rho_n L\left(   L_2^2 \tau^2 + 2 \frac{\Delta^2}{\tau^2}\right) \nonumber\\&& + 8 n^{1/q + 1/2} \sqrt{\rho_n}LD_p\left(L\tau + \frac{2\Delta}{\tau}\right);
    \end{eqnarray}
    \item under Assumptions 1, 2(s), 3 and with $\gamma \leq \frac{\mu}{96 n^{\nicefrac{2}{q}}\rho_nL^2}$:
    \begin{eqnarray}
        \label{teor1.2_a}
        \EE\left[V_{z_{N+1}}(z^*)\right]&\leq& V_{z_0}(z^*) \exp\left(-\frac{\mu^2 N}{400n^{\nicefrac{2}{q}}\rho_n L^2}\right) + \nonumber\\&&+ \frac{24n^{2/q} \rho_n}{\mu^2 N}\left(\|F(z^*)\|^2_2 + \sigma^2 \right)  \\&&+ \frac{4n^{2/q + 1} \rho_n }{\mu^2 N}\left(L_2^2 \tau^2 + 2 \frac{ \Delta^2}{\tau^2}\right)\\&& + \frac{4 n^{1/q + 1/2} \sqrt{\rho_n}  D_p}{\gamma \mu^2 N}\left(L\tau + \frac{2\Delta}{\tau}\right).
    \end{eqnarray}
\end{itemize}
\EndProof

\textbf{Proof of \eqref{teor1.1_a}.}
We begin with descent lemma \eqref{DL_a}:
\begin{eqnarray*}
    \gamma\langle g_d(z_k, e_k, \tau, \xi_k), z_k - u \rangle &\leq&  V_{z_k}(u) - V_{z_{k+1}}(u) + \frac{\gamma^2}{2}\|g_d(z_k, e_k, \tau, \xi_k)\|_q^2.
\end{eqnarray*}
Taking $u = z^*$ and using convexity - concavity of $f(x,y)$ in form $\langle F(z^*), z_k - z^* \rangle \geq 0$, we get
\begin{eqnarray*}
    \gamma\langle F(z_k) - F(z^*), z_k &-& u \rangle \leq  V_{z_k}(z^*) - V_{z_{k+1}}(z^*)\\ &+& \gamma\langle F(z_k) -  g_d(z_k, e_k, \tau, \xi_k), z_k - z^* \rangle + \frac{\gamma^2}{2}\|g_d(z_k, e_k, \tau, \xi_k)\|_q^2.
\end{eqnarray*}
With \eqref{SP_fc}, this gives
\begin{eqnarray*}
    \frac{\gamma}{L} \|F(z_k) - F(z^*)\|^2_2 &\leq&  V_{z_k}(z^*) - V_{z_{k+1}}(z^*)\\ &&+ \gamma\langle F(z_k) -  g_d(z_k, e_k, \tau, \xi_k), z_k - u \rangle + \frac{\gamma^2}{2}\|g_d(z_k, e_k, \tau, \xi_k)\|_q^2.
\end{eqnarray*}
Taking full expectation and using Hölder's inequality, \eqref{lemma1_2_a}, \eqref{lemma1_2_2_a}, we have
\begin{eqnarray*}
    \frac{\gamma}{L} \EE\left[\|F(z_k) - F(z^*)\|^2_2\right] &\leq&  \EE\left[V_{z_k}(z^*)\right] - \EE\left[V_{z_{k+1}}(u)\right] \\
    &&+ 2\gamma\left(2 n^{1/q + 1/2} \sqrt{\rho_n} L\tau + 4 n^{1/q + 1/2} \sqrt{\rho_n} \frac{\Delta}{\tau}\right)D_p\\&&+ \frac{\gamma^2}{2}\left(48 n^{2/q} \rho_n \EE\left[\|F(z_k) - F(z^*)\|^2_2\right] + 48 n^{2/q} \rho_n \|F(z^*)\|^2_2\right) \\&&+ \frac{\gamma^2}{2}\left(48n^{2/q} \rho_n \sigma^2 + 8n^{2/q + 1} \rho_n  L_2^2 \tau^2 + 16\frac{ n^{2/q + 1} \rho_n \Delta^2}{\tau^2}\right).
\end{eqnarray*}
$\gamma \leq \nicefrac{1}{48n^{\frac{2}{q}}\rho_n L}$ gives
\begin{eqnarray*}
    \frac{\gamma}{2L}\EE\left[\|F(z_k) - F(z^*)\|^2_2\right] &\leq&  \EE\left[V_{z_k}(z^*)\right] - \EE\left[V_{z_{k+1}}(z^*)\right] \\
    &&+ 2\gamma\left(2 n^{1/q + 1/2} \sqrt{\rho_n} L\tau + 4 n^{1/q + 1/2} \sqrt{\rho_n} \frac{\Delta}{\tau}\right)D_p\\&&+ \frac{\gamma^2}{2}\left(48 n^{2/q} \rho_n \|F(z^*)\|^2_2 + 48 n^{2/q} \rho_n \sigma^2 \right)\\&&+ \frac{\gamma^2}{2}\left( 8n^{2/q + 1} \rho_n  L_2^2 \tau^2 + 16 \frac{ n^{2/q + 1} \rho_n \Delta^2}{\tau^2}\right).
\end{eqnarray*}
It remains to sum up from $k = 0$ to $k = N$:
\begin{eqnarray*}
    \frac{1}{N+1}\sum^N_{k = 0}\EE\left[\|F(z_k) - F(z^*)\|^2_2\right] &\leq&\frac{2LD^2_p}{\gamma N} + 48 \gamma n^{2/q} \rho_nL\left(\|F(z^*)\|^2_2 + \sigma^2\right) \nonumber\\&& + 8\gamma n^{2/q + 1} \rho_n L\left(   L_2^2 \tau^2 + 2 \frac{\Delta^2}{\tau^2}\right) \nonumber\\&& + 8 n^{1/q + 1/2} \sqrt{\rho_n}LD_p\left(L\tau + \frac{2\Delta}{\tau}\right).
\end{eqnarray*}
\EndProof

\textbf{Proof of \eqref{teor1.2_a}.}
Similarly to the previous proof, we begin with descent lemma \eqref{DL_a}:
\begin{eqnarray*}
    \gamma\langle g(z_k, e_k, \tau, \xi_k), z_k - u \rangle &\leq&  V_{z_k}(u) - V_{z_{k+1}}(u) + \frac{\gamma^2}{2}\|g_d(z_k, e_k, \tau, \xi_k)\|_q^2.
\end{eqnarray*}
Taking $u = z^*$ and using $\langle F(z^*), z_k - z^* \rangle \geq 0$, we get:
\begin{eqnarray*}
    \gamma\langle F(z_k) - F(z^*), z_k &-& z^* \rangle \leq  V_{z_k}(z^*) - V_{z_{k+1}}(z^*)\\ &+& \gamma\langle F(z_k) -  g_d(z_k, e_k, \tau, \xi_k), z_k - u \rangle + \frac{\gamma^2}{2}\|g(z_k, e_k, \tau, \xi_k)\|_q^2.
\end{eqnarray*}
With \eqref{lemma2_a}, it gives
\begin{eqnarray*}
    \frac{\gamma\mu}{2}V_{z_k}(z^*)&\leq&  V_{z_k}(z^*) - V_{z_{k+1}}(z^*)\\ &&+ \gamma\langle F(z_k) -  g_d(z_k, e_k, \tau, \xi_k), z_k - u \rangle + \frac{\gamma^2}{2}\|g_d(z_k, e_k, \tau, \xi_k)\|_q^2.
\end{eqnarray*}
Taking full expectation and using  \eqref{lemma1_2_a}, \eqref{lemma1_2_2_a}, we have
\begin{eqnarray*}
   \EE\left[V_{z_{k+1}}(z^*)\right]&\leq& \left(1 - \frac{\gamma\mu}{2}\right)\EE\left[V_{z_k}(z^*)\right]  + 2\gamma\left(2 n^{1/q + 1/2} \sqrt{\rho_n} L\tau + 4 n^{1/q + 1/2} \sqrt{\rho_n} \frac{\Delta}{\tau}\right)D_p\\&&+ \frac{\gamma^2}{2}\left(48 n^{2/q} \rho_n \EE\left[\|F(z_k) - F(z^*)\|^2_2\right] + 48 n^{2/q} \rho_n \|F(z^*)\|^2_2\right) \\&&+ \frac{\gamma^2}{2}\left(48 n^{2/q} \rho_n \sigma^2 + 8n^{2/q + 1} \rho_n  L_2^2 \tau^2 + 16 \frac{ n^{2/q + 1} \rho_n \Delta^2}{\tau^2}\right).
\end{eqnarray*}
Using \eqref{SP_c} and assuming $\gamma \leq \nicefrac{\mu}{\left(96n^{\nicefrac{2}{q}}\rho_n L^2\right)}$:
\begin{eqnarray*}
   \EE\left[V_{z_{k+1}}(z^*)\right]&\leq& \left(1 - \frac{\gamma\mu}{4}\right)\EE\left[V_{z_k}(z^*)\right]  + 2\gamma^2\left(\frac{2 n^{1/q + 1/2} \sqrt{\rho_n} L\tau}{\gamma} +  \frac{4 n^{1/q + 1/2} \sqrt{\rho_n} \Delta}{\gamma \tau} \right)D_p\\&&+ \gamma^2\left(24 n^{2/q} \rho_n \|F(z^*)\|^2_2 + 24 n^{2/q} \rho_n \sigma^2 \right) \\&&+ \gamma^2\left(4n^{2/q + 1} \rho_n  L_2^2 \tau^2 + 8 \frac{ n^{2/q + 1} \rho_n \Delta^2}{\tau^2}\right).
\end{eqnarray*}
It remains to use \eqref{stich} and get
\begin{eqnarray*}
   \EE\left[V_{z_{N+1}}(z^*)\right]&\leq& V_{z_0}(z^*) \exp\left(-\frac{\mu^2 N}{400n^{\nicefrac{2}{q}}\rho_n L^2}\right) + \nonumber\\&&+ \frac{24n^{2/q} \rho_n}{\mu^2 N}\left(\|F(z^*)\|^2_2 + \sigma^2 \right)  \\&&+ \frac{4n^{2/q + 1} \rho_n }{\mu^2 N}\left(L_2^2 \tau^2 + 2 \frac{ \Delta^2}{\tau^2}\right)\\&& + \frac{4 n^{1/q + 1/2} \sqrt{\rho_n}  D_p}{\gamma \mu^2 N}\left(L\tau + \frac{2\Delta}{\tau}\right).
\end{eqnarray*}
\EndProof

\section{Proof of Theorem 2} \label{add_theorem_2} 

\begin{lemma}\label{lemmt}
Let $z, g, g_{1/2} \in \RR^n$ and $\mathcal{Z} \subset \RR^n$. Then for $z_{1/2} = \text{prox}_z(g)$ and $z_{1} = \text{prox}_z(g_{1/2})$ and for all $u \in \mathcal{Z}$ we have
\begin{equation}
    \label{tlemma}
    \langle g_{1/2}, z_{1/2} - u\rangle  \leq V_{z}(u) - V_{z_1}(u) +  \frac{1}{2}\|g - g_{1/2}\|^2_q - V_{z}(z_{1/2}).
\end{equation}
\end{lemma}

\textbf{Proof.} Using \eqref{lemma3_1} with $z = z$, $z_+ = z_1$, $w = g_{1/2}$, $u = u$ and with $z = z$, $z_+ = z_{1/2}$, $w = g$, $u = z_1$:
\begin{eqnarray*}
         \langle g_{1/2}, z_{1} - u\rangle &\leq& V_{z}(u) - V_{z_1}(u) - V_{z}(z_1), \\
         \langle g, z_{1/2} - z_{1}\rangle &\leq& V_{z}(z_{1}) - V_{z_{1/2}}(z_1) - V_{z}(z_{1/2}).
\end{eqnarray*}
By summing these two inequalities, we get
\begin{eqnarray*}
         \langle g_{1/2}, z_{1/2} - u\rangle  &\leq& V_{z}(u) - V_{z_1}(u) +  \langle g - g_{1/2}, z_{1} - z_{1/2}\rangle \\&&- V_{z_{1/2}}(z_1) - V_{z}(z_{1/2}).
\end{eqnarray*}
Applying Cauchy-Schwartz inequality and property: $V_{z_{1/2}}(z_1) \geq \nicefrac{1}{2} \|z_{1/2} - z_1 \|^2$, we have
\begin{eqnarray*}
         \langle g_{1/2}, z_{1/2} - u\rangle  &\leq& V_{z}(u) - V_{z_1}(u) +  \frac{1}{2}\|g - g_{1/2}\|^2_q - V_{z}(z_{1/2}).
\end{eqnarray*}
\EndProof

\textbf{Theorem.}
\begin{itemize}
    \item By Algorithm 2  with Full coordinates oracle under Assumptions 1, 2, 3 and with $\gamma \leq \nicefrac{1}{2L}$, we have
    \begin{eqnarray}
    \label{teor2.1_a}
        \EE\left[\varepsilon_{sad}(\bar z_N) \right] &\leq& \frac{2D_p^2}{\gamma (N+1)} + 11\gamma\left(n L_2^2 \tau^2 + \sigma^2 + 2 \frac{n\Delta^2}{\tau^2}\right) \nonumber\\
        &&+  2D_p\left( \sqrt{n} L \tau + \frac{2\sqrt{n}\Delta}{\tau}\right),
    \end{eqnarray}
    where 
\begin{equation*}
    \varepsilon_{sad}(\bar z_N) = \max_{y' \in \mathcal{Y}} f(\bar x_N, y') - \min_{x' \in \mathcal{X}} f(x', \bar y_N), 
\end{equation*}
$\bar x_N$, $\bar y_N$ are defined the same way as $\bar z_N$.
    \item By Algorithm 3 with Full coordinates oracle under Assumptions 1, 2(s), 3 and with $p=2$ ($V_x(y) = \nicefrac{1}{2}\|x-y\|^2_2$), $\gamma \leq \nicefrac{1}{6L}$:
    \begin{eqnarray}
        \label{teor2.2_a}
        \EE\big[\|z_{N+1} - z^*\|^2_2\big] &\leq& \exp\left(- \frac{\mu N}{12 L}\right) \left(\|z_0 - z^*\|^2_2 + \|g_f(z_{0},\tau, \xi_{0}) - g_f(z_{0},\tau, \xi_{0})\|^2_2\right)\nonumber\\
    &&+ \frac{1}{\mu^2 N} 12\left(\sigma^2 + nL^2_2\tau^2+ \frac{2n\Delta^2}{\tau^2}\right) \nonumber\\
    &&+\frac{1}{\mu^2 N}\frac{4D_2}{\gamma} \left(\sqrt{n}L\tau+\frac{2\sqrt{n}\Delta}{\tau}\right).
    \end{eqnarray}
\end{itemize}

\textbf{Proof of \eqref{teor2.1_a}.}
We begin with \eqref{tlemma} and taking $z = z_k$, $g = \gamma g_f(z_k, e_k, \tau, \xi_k)$, $g_{1/2} = \gamma g_f(z_{k+1/2}, e_{k+1/2}, \tau, \xi_{k+1/2})$, then $z_{1/2} = z_{k+1/2}$, $z_{1} = z_{k+1}$ and have

\begin{eqnarray*}
    \gamma \langle g_f(z_{k+1/2}, e_{k+1/2}, \tau, &\xi_{k+1/2}&), z_{k+1/2} - u \rangle  \\
     &\leq&  V_{z_k}(u) - V_{z_{k+1}}(u) -  V_{z_k}(z_{k+1/2}) \nonumber \\
     &&+ \frac{\gamma^2}{2}\|g_f(z_{k+1/2}, e_{k+1/2}, \tau, \xi_{k+1/2}) - g_f(z_k, e_k, \tau, \xi_k)\|^2_q \nonumber\\
      &\overset{\eqref{eq:squared_sum}}{\leq}& V_{z_k}(u) - V_{z_{k+1}}(u) -  V_{z_k}(z_{k+1/2}) \nonumber \\ 
      &&+ \frac{3\gamma^2}{2}\|F(z_{k+1/2}) - F(z_k)\|^2_q \nonumber \\ 
      &&+ \frac{3\gamma^2}{2}\|g_f(z_{k+1/2}, e_{k+1/2}, \tau, \xi_{k+1/2}) - F(z_{k+1/2})\|^2_q \nonumber \\ 
     &&+ \frac{3\gamma^2}{2}\|g_f(z_k, e_k, \tau, \xi_k) - F(z_k)\|^2_q \nonumber\\
    &\overset{\eqref{SP_c}}{\leq}& V_{z_k}(u) - V_{z_{k+1}}(u) -  V_{z_k}(z_{k+1/2}) \nonumber \\ 
     &&+ \frac{3\gamma^2L^2}{2}\|z_{k+1/2} - z_k\|^2_2 \nonumber \\ 
     &&+ \frac{3\gamma^2}{2}\|g_f(z_{k+1/2}, e_{k+1/2}, \tau, \xi_{k+1/2}) - F(z_{k+1/2})\|^2_q \nonumber \\ 
     &&+ \frac{3\gamma^2}{2}\|g_f(z_k, e_k, \tau, \xi_k) - F(z_k)\|^2_q.
\end{eqnarray*}
Applying the property: $V_{z_{k}}(z_{k+1/2}) \geq \nicefrac{1}{2} \|z_{k+1/2} - z_k \|^2 \geq \nicefrac{1}{2} \|z_{k+1/2} - z_k \|^2_2$, with $\gamma \leq \nicefrac{1}{2L}$, we get
\begin{eqnarray*}
\gamma \langle g_f(z_{k+1/2}, e_{k+1/2}, &\tau,& \xi_{k+1/2}), z_{k+1/2} - u \rangle  \nonumber \\ &\leq&  V_{z_k}(u) - V_{z_{k+1}}(u) \nonumber \\ &&+ \frac{3\gamma^2}{2}\|g_f(z_{k+1/2}, e_{k+1/2}, \tau, \xi_{k+1/2}) - F(z_{k+1/2})\|^2_q \nonumber \\ &&+ \frac{3\gamma^2}{2}\|g_f(z_k, e_k, \tau, \xi_k) - F(z_k)\|^2_q,
\end{eqnarray*}
and 
\begin{eqnarray*}
\gamma \langle F(z_{k+1/2}), z_{k+1/2} - u \rangle  &\leq&  V_{z_k}(u) - V_{z_{k+1}}(u) \nonumber \\ 
&&+ \gamma \langle F(z_{k+1/2}) - g_f(z_{k+1/2}, e_{k+1/2}, \tau, \xi_{k+1/2}), z_{k+1/2} - u \rangle \nonumber \\ 
&&+ \frac{3\gamma^2}{2}\|g_f(z_{k+1/2}, e_{k+1/2}, \tau, \xi_{k+1/2}) - F(z_{k+1/2})\|^2_q \nonumber \\ 
     &&+ \frac{3\gamma^2}{2}\|g_f(z_k, e_k, \tau, \xi_k) - F(z_k)\|^2_q.
\end{eqnarray*}
Summing over all $k$ from $0$ to $N$, one can have
\begin{eqnarray}
\label{temp404}
\sum\limits_{k=0}^N \langle F(z_{k+1/2}), z_{k+1/2} - u \rangle  &\leq&  \frac{V_{z_0}(u) - V_{z_{K+1}}(u)}{\gamma} \nonumber \\ 
&&+ \sum\limits_{k=0}^N  \langle F(z_{k+1/2}) - g_f(z_{k+1/2}, e_{k+1/2}, \tau, \xi_{k+1/2}), z_{k+1/2} - u \rangle \nonumber \\ 
&&+ \frac{3\gamma}{2} \sum\limits_{k=0}^N \|g_f(z_{k+1/2}, e_{k+1/2}, \tau, \xi_{k+1/2}) - F(z_{k+1/2})\|^2_q \nonumber \\ 
     &&+ \frac{3\gamma}{2} \sum\limits_{k=0}^N \|g_f(z_k, e_k, \tau, \xi_k) - F(z_k)\|^2_q \nonumber \\ 
     &\leq&  \frac{D_p^2}{\gamma} + \sum\limits_{k=0}^N  \langle F(z_{k+1/2}) - g_f(z_{k+1/2}, e_{k+1/2}, \tau, \xi_{k+1/2}), z_{k+1/2} - u \rangle \nonumber \\ 
&&+ \frac{3\gamma}{2} \sum\limits_{k=0}^N \|g_f(z_{k+1/2}, e_{k+1/2}, \tau, \xi_{k+1/2}) - F(z_{k+1/2})\|^2_q \nonumber \\ 
     &&+ \frac{3\gamma}{2} \sum\limits_{k=0}^N \|g_f(z_k, e_k, \tau, \xi_k) - F(z_k)\|^2_q.
\end{eqnarray}

Next we need to connect $\sum\limits_{k=0}^N\langle  F(z_{k+1/2}), z_{k+1/2} - u \rangle$ and $\varepsilon_{sad}(\bar z_{N+1})$. By the definition of $\bar x_N$ and $\bar y_N$, Jensen's inequality and convexity-concavity of $f$:
\begin{eqnarray*}
    \varepsilon_{sad}(\bar z_{N+1})
    &\leq& \max\limits_{y' \in \mathcal{Y}} f\left(\frac{1}{N+1} \left(\sum^N_{k = 0} x_{k+1/2} \right), y'\right) - \min\limits_{x' \in \mathcal{X}} f\left(x', \frac{1}{N+1} \left(\sum^N_{k = 0}  y_{k+1/2} \right)\right) 
    \nonumber \\
    &\leq& \max\limits_{y' \in \mathcal{Y}} \frac{1}{N+1} \sum^N_{k = 0} f(x_{k+1/2}, y')  - \min\limits_{x' \in \mathcal{X}} \frac{1}{N+1} \sum^N_{k = 0} f(x', y_{k+1/2}).
\end{eqnarray*}
Given the fact of linear independence of $x'$ and $y'$:
\begin{eqnarray*}
    \varepsilon_{sad}(\bar z_N) &\leq& \max\limits_{(x', y') \in \mathcal{Z}}\frac{1}{N+1} \sum^N_{k = 0}\left(f(x_{k+1/2}, y')  - f(x', y_{k+1/2}) \right) .
\end{eqnarray*}
Using convexity and concavity of the function $f$:
\begin{eqnarray}
\label{connect_e_sad}
    \varepsilon_{sad}(\bar z_N) &\leq&  \max\limits_{(x', y') \in \mathcal{Z}}\frac{1}{N+1} \sum^N_{k = 0} \left(f(x_{k+1/2}, y')  - f(x', y_{k+1/2}) \right)   \nonumber \\
    &= & \max\limits_{(x', y') \in \mathcal{Z}} \frac{1}{N+1} \sum^N_{k = 0}  \left(f(x_{k+1/2}, y') - f(x_{k+1/2}, y_{k+1/2}) + f(x_{k+1/2}, y_{k+1/2}) - f(x', y_{k+1/2}) \right) \nonumber \\
    &\leq& \max\limits_{(x', y') \in \mathcal{Z}} \frac{1}{N+1} \sum^N_{k = 0}  \left(\langle \nabla_y f (x_{k+1/2}, y_{k+1/2}), y'-y_k \rangle + \langle \nabla_x f (x_{k+1/2}, y_{k+1/2}), x_k-x' \rangle \right) \nonumber \\
    &\leq& \max\limits_{u \in \mathcal{Z}} \frac{1}{N+1} \sum^N_{k = 0}  \langle F(z_{k+1/2}), z_{k+1/2} - u\rangle.
\end{eqnarray}
\eqref{connect_e_sad} together with \eqref{temp404} gives
\begin{eqnarray*}
\varepsilon_{sad}(\bar z_N) &\leq& \frac{D_p^2}{\gamma (N+1)} + \frac{1}{N+1} \max\limits_{u \in \mathcal{Z}} \left[\sum\limits_{k=0}^N  \langle F(z_{k+1/2}) - g_f(z_{k+1/2}, e_{k+1/2}, \tau, \xi_{k+1/2}), z_{k+1/2} - u \rangle\right] \nonumber \\ 
&&+ \frac{3\gamma}{2(N+1)} \sum\limits_{k=0}^N \|g_f(z_{k+1/2}, e_{k+1/2}, \tau, \xi_{k+1/2}) - F(z_{k+1/2})\|^2_q \nonumber \\ 
     &&+ \frac{3\gamma}{2(N+1)} \sum\limits_{k=0}^N \|g_f(z_k, e_k, \tau, \xi_k) - F(z_k)\|^2_q.
\end{eqnarray*}

Taking the full expectation and using \eqref{lemma1_3_a} with \eqref{b}:
\begin{eqnarray}
    \label{344}
    \EE\left[\varepsilon_{sad}(\bar z_N) \right] &\leq& \frac{D_p^2}{\gamma (N+1)} + \frac{1}{N+1} \EE\left[\max\limits_{u \in \mathcal{Z}} \left[\sum\limits_{k=0}^N  \langle F(z_{k+1/2}) - g_f(z_{k+1/2}, e_{k+1/2}, \tau, \xi_{k+1/2}), z_{k+1/2} - u \rangle\right]\right] \nonumber \\ 
&&+ 9\gamma\left(n L_2^2 \tau^2 + \sigma^2 + 2 \frac{n\Delta^2}{\tau^2}\right).
\end{eqnarray}

To finish the proof it remains to estimate \\
$\EE\left[\max\limits_{u \in \mathcal{Z}} \left[\sum\limits_{k=0}^N  \langle F(z_{k+1/2}) - g_f(z_{k+1/2}, e_{k+1/2}, \tau, \xi_{k+1/2}), z_{k+1/2} - u \rangle\right]\right]$. Let define sequence $v$: $v_0 \eqdef z_{1/2}$, $v_{k+1} \eqdef \text{prox}_{v_k} (-\gamma \delta_k)$ with \\ $\delta_k = g_f(z_{k+1/2}, e_{k+1/2}, \tau, \xi_{k+1/2})-F(z_{k+1/2}) $:
\begin{eqnarray}
\label{temp177}
    \sum\limits_{k=0}^N \langle -\delta_k, z_{k+1/2} - u \rangle = \sum\limits_{k=0}^N \langle -\delta_k, z_{k+1/2} - v_k \rangle + 
    \sum\limits_{k=0}^N \langle -\delta_k,  v_k - u \rangle . 
\end{eqnarray}
By the definition of $v$ and an optimal condition for the prox-operator, we have for all $u \in \mathcal{Z}$
\begin{eqnarray*}
    \langle -\gamma\delta_k - \nabla d(v_{k+1})  + \nabla d(v_{k+1}), u - v_{k+1} \rangle \geq 0.
\end{eqnarray*}
Rewriting this inequality, we get
\begin{eqnarray*}
    \langle -\gamma\delta_k, v_k  - u \rangle &\leq& \langle -\gamma\delta_k, v_k - v_{k+1} \rangle  + \langle \nabla d (v_{k+1}) - \nabla d (v_k), u - v_{k+1} \rangle \nonumber\\
&\leq& \langle -\gamma\delta_k, v_k - v_{k+1} \rangle + V_{v_k}(u) -  V_{v_{k+1}}(u) - V_{v_k}(v_{k+1}).
\end{eqnarray*}
Bearing in mind the Bregman divergence property $2V_x(y) \geq \|x-y\|_p^2$:
\begin{eqnarray*}
    \langle -\gamma\delta_k, v_k  - u \rangle \leq \langle -\gamma\delta_k, v_k - v_{k+1} \rangle + V_{v_k}(u) -  V_{v_{k+1}}(u) - \frac{1}{2}\|v_{k+1} - v_k\|_p^2.
\end{eqnarray*}
Using the definition of the conjugate norm:
\begin{eqnarray*}
    \langle -\gamma\delta_k, v_k  - u \rangle &\leq& \|\gamma \delta_k\|_q\cdot\| v_k - v_{k+1} \|_p + V_{v_k}(u) -  V_{v_{k+1}}(u) - \frac{1}{2}\|v_{k+1} - v_k\|_p^2 \nonumber\\ 
    &\leq& \frac{\gamma^2}{2}\|\delta_k\|_q^2 + V_{v_k}(u) -  V_{v_{k+1}}(u).
\end{eqnarray*}
Summing over $k$ from $0$ to $N$:
\begin{eqnarray}
    \gamma \sum\limits_{k=0}^N  \langle -\delta_k,  v_k - u \rangle &\leq& V_{v_0}(u) -  V_{v_{N+1}}(u) + \frac{\gamma^2}{2}\sum\limits_{k=0}^N \|\delta_k\|_q^2 \nonumber\\
\label{temp192}
    &\leq& D_p^2 + \frac{\gamma^2}{2}\sum\limits_{k=0}^N \|\delta_k\|_q^2.
\end{eqnarray}
Substituting \eqref{temp192} into \eqref{temp177}:
\begin{eqnarray*}
     \sum\limits_{k=0}^N \langle -\delta_k, z_{k+1/2} - u \rangle = \sum\limits_{k=0}^N \langle \delta_k, v_k - z_{k+1/2}  \rangle +
    \frac{D_p^2}{\gamma} + \frac{\gamma}{2}\sum\limits_{k=0}^N \|\delta_k\|_q^2.  
\end{eqnarray*}
The right side is independent of $u$, then
\begin{eqnarray*}
     \max_{u \in \mathcal{Z}} \sum\limits_{k=0}^N \langle -\delta_k, z_{k+1/2} - u \rangle &\leq& \sum\limits_{k=0}^N \langle \delta_k, v_k - z_{k+1/2}  \rangle + 
    \frac{D_p^2}{\gamma} + \frac{\gamma^2}{2}\sum\limits_{k=0}^N \|\delta_k\|_q^2.  
\end{eqnarray*}
Taking the full expectation with independence $v_k - z_{k+1/2} $, $\xi_{k+1/2}, e_{k+1/2}$  and using \eqref{lemma1_3_a}, \eqref{lemma1_3_1_a}, we get
\begin{eqnarray}
    \label{405}
     \EE\left[\max_{u \in \mathcal{Z}} \sum\limits_{k=0}^N \langle -\delta_k, z_{k+1/2} - u \rangle\right] &\leq& \EE\left[\sum\limits_{k=0}^N \langle \delta_k, v_k - z_{k+1/2}  \rangle\right] + 
    \frac{D_p^2}{\gamma} + \frac{\gamma}{2}\sum\limits_{k=0}^N \EE\left[\|\delta_k\|_q^2\right] \nonumber\\
    &\leq& \EE\left[\sum\limits_{k=0}^N \langle \EE_{e_{k+1/2}, \xi_{k+1/2}}\left[\delta_k\right], v_k - z_{k+1/2}  \rangle\right] + 
    \frac{D_p^2}{\gamma} + \frac{\gamma}{2}\sum\limits_{k=0}^N \EE\left[\|\delta_k\|_q^2\right] \nonumber\\
    &\leq& 2(N+1)D_p\left( \sqrt{n} L \tau + \frac{2\sqrt{n}\Delta}{\tau} \right)+ 
    \frac{D_p^2}{\gamma} \nonumber\\
    &&+ \frac{3\gamma (N+1)}{2} \left( \sigma^2 + nL_2^2 \tau^2 + \frac{2n\Delta^2}{\tau^2}\right).  
\end{eqnarray}
Connecting \eqref{344} and \eqref{405}, we have
\begin{eqnarray*}
    \EE\left[\varepsilon_{sad}(\bar z_N) \right] &\leq& \frac{2D_p^2}{\gamma (N+1)} + 11\gamma\left(n L_2^2 \tau^2 + \sigma^2 + 2 \frac{n\Delta^2}{\tau^2}\right) \nonumber\\
    &&+  2D_p\left( \sqrt{n} L \tau + \frac{2\sqrt{n}\Delta}{\tau}\right).
\end{eqnarray*}

\textbf{Proof of \eqref{teor2.2_a}.}
Similarly to the previous proof, let begin with \eqref{tlemma} and take full expectation:
\begin{eqnarray}
\label{temppp}
    \EE\big[\|z_{k+1} - z^*\|^2_2\big] &\leq& \EE\big[\|z_k - z^*\|^2_2\big] - 2\gamma\EE\big[\langle g_f(z_{k+ \nicefrac{1}{2}},\tau, \xi_{k+ \nicefrac{1}{2}}), z_{k+ \nicefrac{1}{2}} - z^*\rangle\big] \nonumber \\&&+ \gamma^2 \EE\big[\|g_f(z_{k+ \nicefrac{1}{2}},\tau, \xi_{k+ \nicefrac{1}{2}}) - g_f(z_{k- \nicefrac{1}{2}},\tau, \xi_{k- \nicefrac{1}{2}})\|^2_2\big]  \nonumber \\&&- \EE\big[\|z_{k+ \nicefrac{1}{2}} - z_{k}\|^2_2\big].
\end{eqnarray}
Next we work with $\EE\left[\|g_f(z_{k+ \nicefrac{1}{2}},\tau, \xi_{k+ \nicefrac{1}{2}}) - g_f(z_{k- \nicefrac{1}{2}},\tau, \xi_{k- \nicefrac{1}{2}})\|^2_2\right]$:
\begin{eqnarray*}
    \EE\big[\|g_f(z_{k+ \nicefrac{1}{2}},\tau, \xi_{k+ \nicefrac{1}{2}}) - g_f&(z_{k- \nicefrac{1}{2}},&\tau, \xi_{k- \nicefrac{1}{2}})\|^2_2\big]\\  &\overset{\eqref{eq:squared_sum}}{\leq}& 3\EE\left[\|g_f(z_{k+ \nicefrac{1}{2}},\tau, \xi_{k+ \nicefrac{1}{2}}) - F(z_{k+ \nicefrac{1}{2}})\|^2_2\right]\\&& + 3\EE\left[\|g_f(z_{k-\nicefrac{1}{2}},\tau, \xi_{k+ \nicefrac{1}{2}}) - F(z_{k- \nicefrac{1}{2}})\|^2_2\right] \\&&+
    3\EE\left[\|F(z_{k+ \nicefrac{1}{2}}) - F(z_{k- \nicefrac{1}{2}})\|^2_2\right]\\ &\overset{\eqref{lemma1_3_a},\eqref{SP_c}}{\leq}& 3L^2\EE\left[\|z_{k+ \nicefrac{1}{2}}- z_{k- \nicefrac{1}{2}}\|^2_2\right] +6\left(\sigma^2 + nL^2_2\tau^2+ \frac{2n\Delta^2}{\tau^2}\right)\\ &\overset{\eqref{eq:squared_sum}}{\leq}& 6L^2\EE\left[\|z_{k+ \nicefrac{1}{2}}- z_{k}\|^2_2\right] + 6L^2\EE\left[\|z_{k}- z_{k- \nicefrac{1}{2}}\|^2_2\right]\\&&+6\left(\sigma^2 + nL^2_2\tau^2+ \frac{2n\Delta^2}{\tau^2}\right) \\&\leq& 6L^2\EE\left[\|z_{k+ \nicefrac{1}{2}}- z_{k}\|^2_2\right] \\&&+ 6\gamma^2L^2\EE\left[\|g_f(z_{k- \nicefrac{1}{2}},\tau, \xi_{k- \nicefrac{1}{2}}) - g_f(z_{k- \nicefrac{3}{2}},\tau, \xi_{k- \nicefrac{3}{2}})\|^2_2\right]\\&&+6\left(\sigma^2 + nL^2_2\tau^2+ \frac{2n\Delta^2}{\tau^2}\right).
\end{eqnarray*}
In last inequality we use non-expansiveness of Euclidean prox operator. By simple transformation:
\begin{eqnarray*}
    \EE\big[\|g_f(z_{k+ \nicefrac{1}{2}},&\tau,& \xi_{k+ \nicefrac{1}{2}}) - g_f(z_{k- \nicefrac{1}{2}},\tau, \xi_{k- \nicefrac{1}{2}})\|^2_2\big] \\&\leq& 12L^2\EE\left[\|z_{k+ \nicefrac{1}{2}}- z_{k}\|^2_2\right] \\&&+ 12\gamma^2L^2\EE\left[\|g_f(z_{k- \nicefrac{1}{2}},\tau, \xi_{k- \nicefrac{1}{2}}) - g_f(z_{k- \nicefrac{3}{2}},\tau, \xi_{k- \nicefrac{3}{2}})\|^2_2\right]\\
    &&-\EE\big[\|g_f(z_{k+ \nicefrac{1}{2}},\tau, \xi_{k+ \nicefrac{1}{2}}) - g_f(z_{k- \nicefrac{1}{2}},\tau, \xi_{k- \nicefrac{1}{2}})\|^2_2\big]\\&&
    +12\left(\sigma^2 + nL^2_2\tau^2+ \frac{2n\Delta^2}{\tau^2}\right).
\end{eqnarray*}
If $\gamma \leq \nicefrac{1}{6L}$, then $12\gamma^2L^2 \leq 1-\mu\gamma$, and we can rewrite previous inequality:
\begin{eqnarray}
\label{temppp1}
    \EE\big[\|g_f(z_{k+ \nicefrac{1}{2}},&\tau,& \xi_{k+ \nicefrac{1}{2}}) - g_f(z_{k- \nicefrac{1}{2}},\tau, \xi_{k- \nicefrac{1}{2}})\|^2_2\big] \nonumber\\&\leq& 12L^2\EE\left[\|z_{k+ \nicefrac{1}{2}}- z_{k}\|^2_2\right] \nonumber\\&&+ (1-\mu\gamma)\EE\left[\|g_f(z_{k- \nicefrac{1}{2}},\tau, \xi_{k- \nicefrac{1}{2}}) - g_f(z_{k- \nicefrac{3}{2}},\tau, \xi_{k- \nicefrac{3}{2}})\|^2_2\right]\nonumber\\
    &&-\EE\big[\|g_f(z_{k+ \nicefrac{1}{2}},\tau, \xi_{k+ \nicefrac{1}{2}}) - g_f(z_{k- \nicefrac{1}{2}},\tau, \xi_{k- \nicefrac{1}{2}})\|^2_2\big]\nonumber\\&&
    +12\left(\sigma^2 + nL^2_2\tau^2+ \frac{2n\Delta^2}{\tau^2}\right).
\end{eqnarray}
Next we consider $- 2\gamma\EE\left[\langle g_f(z_{k+ \nicefrac{1}{2}},\tau, \xi_{k+ \nicefrac{1}{2}}), z_{k+ \nicefrac{1}{2}} - z^*\rangle\right]$:
\begin{eqnarray}
\label{temppp2}
    - 2\gamma\EE\big[\langle g_f(z_{k+ \nicefrac{1}{2}},&\tau,& \xi_{k+ \nicefrac{1}{2}}), z_{k+ \nicefrac{1}{2}} - z^*\rangle\big] \nonumber \\
    &=& - 2\gamma\EE\left[\langle F(z_{k+ \nicefrac{1}{2}}), z_{k+ \nicefrac{1}{2}} - z^*\rangle\right] \nonumber\\&&+2\gamma\EE\left[\langle F(z_{k+ \nicefrac{1}{2}}) -  g_f(z_{k+ \nicefrac{1}{2}},\tau, \xi_{k+ \nicefrac{1}{2}}), z_{k+ \nicefrac{1}{2}} - z^*\rangle\right]\nonumber\\&\leq&  - 2\gamma\EE\left[\langle F(z_{k+ \nicefrac{1}{2}}), z_{k+ \nicefrac{1}{2}} - z^*\rangle\right] \nonumber\\&&+4\gamma\|\EE \left[ F(z_{k+ \nicefrac{1}{2}}) - g_f(z_{k+ \nicefrac{1}{2}},\tau, \xi_{k+ \nicefrac{1}{2}})\right]\|_2 D_2\nonumber\\&\overset{\eqref{lemma1_3_1_a}}{\leq}& -2\gamma\EE\left[\langle F(z_{k+ \nicefrac{1}{2}}), z_{k+ \nicefrac{1}{2}} - z^*\rangle\right] \nonumber\\&&+4\gamma\left(\sqrt{n}L\tau+\frac{2\sqrt{n}\Delta}{\tau}\right)D_2 \nonumber\\
    &\overset{\eqref{lemma2_a}}{\leq}& -2\gamma \mu\EE\left[\|z_{k+ \nicefrac{1}{2}} - z^*\|^2_2\right] \nonumber \\&&+4\gamma\left(\sqrt{n}L\tau+\frac{2\sqrt{n}\Delta}{\tau}\right)D_2 \nonumber\\
    &\leq& -\gamma \mu\EE\left[\|z_{k} - z^*\|^2_2\right] + 2 \gamma \mu\EE\left[\|z_{k+ \nicefrac{1}{2}} - z_k\|^2_2\right] \nonumber \\&&+4\gamma\left(\sqrt{n}L\tau+\frac{2\sqrt{n}\Delta}{\tau}\right)D_2.
\end{eqnarray}
Combining \eqref{temppp}, \eqref{temppp1}, and \eqref{temppp2}, we have 
\begin{eqnarray*}
    \EE\big[\|z_{k+1} - z^*\|^2_2\big] &+&\EE\big[\|g_f(z_{k+ \nicefrac{1}{2}},\tau, \xi_{k+ \nicefrac{1}{2}}) - g_f(z_{k- \nicefrac{1}{2}},\tau, \xi_{k- \nicefrac{1}{2}})\|^2_2\big] \\
    &\leq& (1 - \gamma \mu)\left(\EE\big[\|z_k - z^*\|^2_2\big] + \EE\left[\|g_f(z_{k- \nicefrac{1}{2}},\tau, \xi_{k- \nicefrac{1}{2}}) - g_f(z_{k- \nicefrac{3}{2}},\tau, \xi_{k- \nicefrac{3}{2}})\|^2_2\right]\right) \\
    &&+ (2 \gamma \mu + 12\gamma^2L^2 - 1)\EE\left[\|z_{k+ \nicefrac{1}{2}} - z_k\|^2_2\right] \nonumber \\&&
    +\gamma^2 \left[12\left(\sigma^2 + nL^2_2\tau^2+ \frac{2n\Delta^2}{\tau^2}\right) +\frac{4D_2}{\gamma} \left(\sqrt{n}L\tau+\frac{2\sqrt{n}\Delta}{\tau}\right)\right].
\end{eqnarray*}
With $\gamma \leq \nicefrac{1}{6L}$ we have $12\gamma^2L^2 \leq 1-2\mu\gamma$ and 
\begin{eqnarray*}
    \EE\big[\|z_{k+1} - z^*\|^2_2\big] &+&\EE\big[\|g_f(z_{k+ \nicefrac{1}{2}},\tau, \xi_{k+ \nicefrac{1}{2}}) - g_f(z_{k- \nicefrac{1}{2}},\tau, \xi_{k- \nicefrac{1}{2}})\|^2_2\big] \\
    &\leq& (1 - \gamma \mu)\left(\EE\big[\|z_k - z^*\|^2_2\big] + \EE\left[\|g_f(z_{k- \nicefrac{1}{2}},\tau, \xi_{k- \nicefrac{1}{2}}) - g_f(z_{k- \nicefrac{3}{2}},\tau, \xi_{k- \nicefrac{3}{2}})\|^2_2\right]\right) \\&&
    +\gamma^2 \left[12\left(\sigma^2 + nL^2_2\tau^2+ \frac{2n\Delta^2}{\tau^2}\right) +\frac{4D_2}{\gamma} \left(\sqrt{n}L\tau+\frac{2\sqrt{n}\Delta}{\tau}\right)\right].
\end{eqnarray*}
It remains to apply \eqref{stich} and then :
\begin{eqnarray*}
    \EE\big[\|z_{N+1} - z^*\|^2_2\big] &\leq& \exp\left(- \frac{\mu N}{12 L}\right) \left(\|z_0 - z^*\|^2_2 + \|g_f(z_{0},\tau, \xi_{0}) - g_f(z_{0},\tau, \xi_{0})\|^2_2\right)\\
    &&+ \frac{1}{\mu^2 N} \left[12\left(\sigma^2 + nL^2_2\tau^2+ \frac{2n\Delta^2}{\tau^2}\right) +\frac{4D_2}{\gamma} \left(\sqrt{n}L\tau+\frac{2\sqrt{n}\Delta}{\tau}\right)\right].
\end{eqnarray*}

\EndProof

\section{Other approach for $e$ in Algorithm \ref{alg2}} \label{same direction}

This algorithm is an easy modification of Algorithm \ref{alg2}. The only difference is that we use the same direction $e$ and random variable $\xi$ within one iteration \\
\begin{minipage}{1.\linewidth}
     \begin{algorithm}[H]
	\caption{{\tt zoESVIA (same direction)}}
	\label{alg2_1}
\begin{algorithmic}
\State 
\noindent {\bf Input:} $z_0$, $N$, $\gamma$, $\tau$.
\State Choose oracle $\text{grad}$ from $G, g_d, g_f.$.
\For {$k=0,1, 2, \ldots, N$ }
    \State Sample indep. $e_k$, $\xi_k$.
    \State $d_{k} =  \text{grad}(z_{k}, e_{k}, \tau, \xi_k)$.
    \State $z_{k+1/2} = \text{prox}_{z_k}(\gamma \cdot d_{k})$.
    \State $d_{k+1/2} =  \text{grad}(z_{k+1/2}, e_{k}, \tau, \xi_{k})$.
    \State $z_{k+1} = \text{prox}_{z_k}(\gamma \cdot d_{k+1/2})$.
\EndFor
\State 
\noindent {\bf Output:} $z_{N+1}$ or $\bar z_{N+1}$.
\end{algorithmic}
\end{algorithm}
\end{minipage}
In this section we consider euclidean setup: $V_x(y) = \nicefrac{1}{2}\|x - y\|^2_2$. Used approach is based on \cite{gidel2018variational}.

\textbf{Theorem.}

By Algorithm 4 with Random direction oracle under Assumptions 1, 2, 3 and $\gamma \leq \nicefrac{1}{2nL_2}$, we get
\begin{eqnarray}
\label{teor4.1_a}
     \EE[\varepsilon_{sad}(\bar z_{N})]&\leq&  \frac{D^2_2}{\gamma N} + 210\gamma n^2 L^2_2D^2_2 \nonumber\\&&+ 24\gamma\left(n^2L^2\tau^2 + \frac{n^2\Delta^2}{\tau^2}\right) +  12\left(n L\tau + \frac{n\Delta}{\tau}\right)D_2 \nonumber\\&&+ 200\gamma\left(n\EE\left[\|F(z^*)\|^2_2\right] + \frac{n\sigma^2}{2}\right),
\end{eqnarray}
    where \begin{equation*}
    \varepsilon_{sad}(\bar z_N) = \max_{y' \in \mathcal{Y}} f(\bar x_N, y') - \min_{x' \in \mathcal{X}} f(x', \bar y_N). 
\end{equation*}

\textbf{Proof of \eqref{teor4.1_a}}
We begin with applying Lemma \ref{lemmt}
\begin{eqnarray*}
    \|z_{k+1} - u\|^2_2 &\leq&  \|z_k - u\|^2_2 - 2\langle\gamma g_d(z_{k+ \nicefrac{1}{2}},e_k, \tau, \xi_k ), z_{k+\nicefrac{1}{2}} - u\rangle  \\&&+ \gamma^2 \|g_d(z_{k+ \nicefrac{1}{2}},e_k, \tau, \xi_k ) - g_d(z_{k},e_k, \tau, \xi_k)\|^2_2 - \|z_{k+\nicefrac{1}{2}} - z_{k}\|^2_2
\end{eqnarray*}
Next, using triangle inequality, we have
\begin{eqnarray*}
    \|z_{k+1} - u\|^2_2 &\leq&  \|z_k - u\|^2_2 - 2\langle\gamma g_d(z_{k+ \nicefrac{1}{2}},e_k, \tau, \xi_k ), z_{k+\nicefrac{1}{2}} - u\rangle  \\&&+ \gamma^2 \|g_d(z_{k+ \nicefrac{1}{2}},e_k, \tau, \xi_k ) -n\langle F(z_{k+ \nicefrac{1}{2}}, \xi_k), e_k \rangle e_k \|^2_2 \\&&+\gamma^2\|g_d(z_{k},e_k, \tau, \xi_k) - n\langle F(z_k, \xi_k), e_k \rangle e_k\|^2_2  \\&&+ \gamma^2\|n\langle F(z_{k+ \nicefrac{1}{2}}, \xi_k), e_k \rangle e_k  - n\langle F(z_k, \xi_k), e_k \rangle e_k\|^2_2 -\|z_{k+\nicefrac{1}{2}} - z_{k}\|^2_2
\end{eqnarray*}
Using \ref{SP_c}, we get
\begin{eqnarray*}
    \|z_{k+1} - u\|^2_2 &\leq&  \|z_k - u\|^2_2 - 2\langle\gamma g_d(z_{k+ \nicefrac{1}{2}},e_k, \tau, \xi_k ), z_{k+\nicefrac{1}{2}} - u\rangle  \\&&+ \gamma^2 \|g_d(z_{k+ \nicefrac{1}{2}},e_k, \tau, \xi_k ) -n\langle F(z_{k+ \nicefrac{1}{2}}, \xi_k), e_k \rangle e_k \|^2_2 \\&&+\gamma^2\|g_d(z_{k},e_k, \tau, \xi_k) - n\langle F(z_k, \xi_k), e_k \rangle e_k\|^2_2  \\&&+ (\gamma^2n^2 L^2(\xi_k) - 1)\|z_{k+ \nicefrac{1}{2}} - z_k\|^2_2 
\end{eqnarray*}
By simple transformation we rewrite previous inequality 
\begin{eqnarray*}
     \langle F(z_{k+\nicefrac{1}{2}}), z_{k+\nicefrac{1}{2}} - u\rangle&\leq&  \|z_k - u\|^2_2 - \|z_{k+1} - u\|^2_2 \\&&- 2\gamma\langle g_d(z_{k+ \nicefrac{1}{2}},e_k, \tau, \xi_k ) - F(z_{k+\nicefrac{1}{2}}, \xi_k), z_{k+\nicefrac{1}{2}} - u\rangle  \\&&+ \gamma^2 \|g_d(z_{k+ \nicefrac{1}{2}},e_k, \tau, \xi_k ) -n\langle F(z_{k+ \nicefrac{1}{2}}, \xi_k), e_k \rangle e_k \|^2_2 \\&&+\gamma^2\|g_d(z_{k},e_k, \tau, \xi_k) - n\langle F(z_k, \xi_k), e_k \rangle e_k\|^2_2  \\&&+ (\gamma^2n^2 L^2_2 - 1)\|z_{k+ \nicefrac{1}{2}} - z_k\|^2_2.
\end{eqnarray*}
We estimate some terms from the right side of the inequality.
\begin{eqnarray*}
     &&\|g_d(z_{k+ \nicefrac{1}{2}},e_k, \tau, \xi_k ) - n\langle F(z_{k+ \nicefrac{1}{2}}, \xi_k), e_k \rangle e_k \|^2_2 \leq \nonumber\\
     &&\frac{n^2}{\tau^2} \left\| \left(
    \begin{array}{c}
    \left( f(x_{k+ \nicefrac{1}{2}} + \tau e_{kx}, y_{k+ \nicefrac{1}{2}},  \xi_k) -   f(x_{k+ \nicefrac{1}{2}}, y_{k+ \nicefrac{1}{2}}, \xi_k) - \langle \nabla_x f(x_{k+ \nicefrac{1}{2}},y_{k+ \nicefrac{1}{2}}, \xi_k), \tau e_{kx}\rangle  \right)e_{kx}\\
    \left( f(x_{k+ \nicefrac{1}{2}}, y_{k+ \nicefrac{1}{2}}, \xi_k) - f(x_{k+ \nicefrac{1}{2}}, y_{k+ \nicefrac{1}{2}} + \tau e_{ky},  \xi_k) +  \langle \nabla_y f(x_{k+ \nicefrac{1}{2}},y_{k+ \nicefrac{1}{2}}, \xi_k), \tau e_{ky}\rangle  \right)e_{ky} 
    \end{array}
    \right)  \right\|^2_2  \\
&&+ \frac{n^2}{\tau^2} \left\| \left(
    \begin{array}{c}
    \left(\delta(x_{k+ \nicefrac{1}{2}} + \tau e_{kx}, y_{k+ \nicefrac{1}{2}}) -  \delta(x_{k+ \nicefrac{1}{2}}, y_{k+ \nicefrac{1}{2}}) \right)e_{kx}\\
    \left(\delta(x_{k+ \nicefrac{1}{2}}, y_{k+ \nicefrac{1}{2}}) - \delta(x_{k+ \nicefrac{1}{2}}, y_{k+ \nicefrac{1}{2}} + \tau e_{ky}) \right)e_{ky} 
    \end{array}
    \right)\right\|^2_2.
\end{eqnarray*}
Using $L$-smoothness of function $f(\cdot)$ and \ref{eq:block_vec}, we note that 
\begin{eqnarray*}
     \|g_d(z_{k+ \nicefrac{1}{2}},e_k, \tau, \xi_k ) - n\langle F(z_{k+ \nicefrac{1}{2}}, \xi_k), e_k \rangle e_k \|^2_2 &\leq& \frac{n^2}{\tau^2} \left(L^2\|\tau e_{kx}\|^2_2 + L^2 \|\tau
     e_{ky}\|^2_2\right) \\&&+ 4\frac{n^2\Delta^2}{\tau^2}\left(\| e_{kx}\|^2_2+  \| e_{ky}\|^2_2\right) \\
     &\leq& 4\left(n^2L^2\tau^2 + \frac{n^2\Delta^2}{\tau^2}\right)
\end{eqnarray*}
Similarly, we estimate the following value. Using $L$-smoothness of function $f(\cdot)$, we have
\begin{eqnarray*}
     \|g_d(z_{k},e_k, \tau, \xi_k ) -n\langle F(z_{k}, \xi_k), e_k \rangle e_k \|^2_2 &\leq&\frac{n^2}{\tau^2} \left(L^2\|\tau e_{kx}\|^2_2 + L\|\tau
     e_{ky}\|^2_2\right) + 4\frac{n^2\Delta^2}{\tau^2}\left(\| e_{kx}\|^2_2+  \| e_{ky}\|^2_2\right) \\
     &\leq& 4\left(n^2L^2\tau^2 + \frac{n^2\Delta^2}{\tau^2}\right)
\end{eqnarray*}

Substituting the previous inequality we have
\begin{eqnarray*}
     \langle F(z_{k+\nicefrac{1}{2}}), z_{k+\nicefrac{1}{2}} - u\rangle&\leq&  \|z_k - u\|^2_2 - \|z_{k+1} - u\|^2_2 \\&&- 2\gamma\langle g_d(z_{k+ \nicefrac{1}{2}},e_k, \tau, \xi_k ) - F(z_{k+\nicefrac{1}{2}}), z_{k+\nicefrac{1}{2}} - u\rangle  \\&&+ 8\gamma^2\left(n^2L^2\tau^2 + \frac{n^2\Delta^2}{\tau^2}\right)  \\&&+ (\gamma^2n^2 L^2_2 - 1)\|z_{k+ \nicefrac{1}{2}} - z_k\|^2_2 
\end{eqnarray*}

Consider $\mathcal{G} = \langle g_d(z_{k+ \nicefrac{1}{2}},e_k, \tau, \xi_k ) - F(z_{k+\nicefrac{1}{2}}),u - z_{k+\nicefrac{1}{2}}\rangle$, by simple transformations we get
\begin{eqnarray*}
     \mathcal{G} &=& \langle g_d(z_{k+ \nicefrac{1}{2}},e_k, \tau, \xi_k ) - g_d(z_{k},e_k, \tau, \xi_k ),u - z_{k+\nicefrac{1}{2}}\rangle + \langle F(z_{k}) -  F(z_{k+\nicefrac{1}{2}}),u - z_{k+\nicefrac{1}{2}}\rangle\\&&+ \langle g_d(z_{k},e_k, \tau, \xi_k ) - F(z_{k}),u - z_k \rangle + \langle g_d(z_{k},e_k, \tau, \xi_k ) - F(z_{k}), z_k - z_{k+\nicefrac{1}{2}}\rangle \\&\leq& 2nL_2\|z_{k} - z_{k+\nicefrac{1}{2}}\|_2\|u - z_{k+\nicefrac{1}{2}}\|_2 + \|g_d(z_{k},e_k, \tau, \xi_k ) - F(z_{k})\|_2\|z_k - z_{k+\nicefrac{1}{2}}\|_2 \\&&+ \|g_d(z_{k+ \nicefrac{1}{2}},e_k, \tau, \xi_k ) -n\langle F(z_{k+ \nicefrac{1}{2}}), e_k \rangle e_k \|_2\|u - z_{k+\nicefrac{1}{2}}\|_2 \\&&+ \|g_d(z_{k},e_k, \tau) - n\langle F(z_k), e_k \rangle e_k\|_2\|u - z_{k+\nicefrac{1}{2}}\|_2 + \langle g_d(z_{k},e_k, \tau, \xi_k ) - F(z_{k}),u - z_k \rangle.
\end{eqnarray*}
Using $2\|a\|\|b\|\leq C\|a\|^2_2 + \frac{1}{C}\|b\|^2_2$
\begin{eqnarray*}
2\gamma \mathcal{G} &\leq& \frac{1}{2}\|z_{k} - z_{k+\nicefrac{1}{2}}\|^2_2 +8\gamma^2 n^2 L^2_2\|u - z_{k+\nicefrac{1}{2}}\|^2_2 + 4\gamma^2\|g_d(z_{k},e_k, \tau, \xi_k ) - F(z_{k})\|^2_2 \\
&&+\frac{1}{4}\|z_k - z_{k+\nicefrac{1}{2}}\|^2_2 \\&&+ 4\gamma\left(n L\tau + \frac{n\Delta}{\tau}\right)\|u - z_{k+\nicefrac{1}{2}}\|_2 + 2\gamma\langle g_d(z_{k},e_k, \tau, \xi_k ) - F(z_{k}),u - z_k \rangle.
\end{eqnarray*}
Summing up we get
\begin{eqnarray*}
     \gamma\langle F(z_{k+\nicefrac{1}{2}}), z_{k+\nicefrac{1}{2}} - u\rangle&\leq&  \|z_k - u\|^2_2 - \|z_{k+1} - u\|^2_2 + 2\gamma\langle g_d(z_{k},e_k, \tau, \xi_k ) - F(z_{k}),u - z_k \rangle \\&&+ 8\gamma^2\left(n^2L^2\tau^2 + \frac{n^2\Delta^2}{\tau^2}\right) +  4\gamma\left(n L\tau + \frac{n\Delta}{\tau}\right)\|u - z_{k+\nicefrac{1}{2}}\|_2\\&&+ (\gamma^2n^2 L^2_2 - \frac{1}{4})\|z_{k+ \nicefrac{1}{2}} - z_k\|^2_2 \\&&+ 8\gamma^2 n^2 L^2_2\|u - z_{k+\nicefrac{1}{2}}\|^2_2 + 4\gamma^2\|g_d(z_{k},e_k, \tau, \xi_k ) - F(z_{k})\|^2_2
\end{eqnarray*}
Assuming $\gamma \leq \nicefrac{1}{2nL_2}$, convexity-concavity of function $f(\dot)$ and summing from $k = 1$ to $k = N$, we get
\begin{eqnarray*}
    \frac{\gamma}{N+1}\sum^N_{k = 0}\langle F(z_{k+\nicefrac{1}{2}}), z_{k+\nicefrac{1}{2}} - u\rangle&\leq&  \frac{D^2_2}{N} + \frac{2\gamma}{N}\sum^N_{k = 1}\langle g_d(z_{k},e_k, \tau, \xi_k ) - F(z_{k}),u - z_k \rangle \\&&+ 8\gamma^2\left(n^2L^2\tau^2 + \frac{n^2\Delta^2}{\tau^2}\right) +  4\gamma\left(n L\tau + \frac{n\Delta}{\tau}\right)D_2 \\&&+ 8\gamma^2 n^2 L^2_2D^2_2 + \frac{4\gamma^2}{N}\sum^N_{k = 1}\|g_d(z_{k},e_k, \tau, \xi_k ) - F(z_{k})\|^2_2
\end{eqnarray*}
Taking full expectation and using \ref{lemma1_2_a} $(q = 2)$, \ref{connect_e_sad} and \ref{SP_c},  we have 
\begin{eqnarray*}
     \gamma\EE[\varepsilon_{sad}(\bar z_{N})]&\leq&  \frac{D^2_2}{N} + 210\gamma^2 n,2 L^2_2D^2_2 \\&&+ 24\gamma^2\left(n^2L^2\tau^2 + \frac{n^2\Delta^2}{\tau^2}\right) +  12\gamma\left(n L\tau + \frac{n\Delta}{\tau}\right)D_2 \\&&+ 200\gamma^2\left(n\EE\left[\|F(z^*)\|^2_2\right] + \frac{n\sigma^2}{2}\right)
\end{eqnarray*}
\EndProof

\section{Proof of Theorem 3}

\begin{theorem} By Algorithm \ref{alg2} under assumption 1, 2, 3 with Mixed oracle $\widetilde{g}_f$ and $\gamma \leq \nicefrac{1}{2L}$, we get
\begin{eqnarray}
\label{teor5.1}
     \EE\left[\varepsilon_{sad}(\bar z_N) \right] &\leq& \frac{D_p^2}{\gamma N} + 2 D_p\left(\sqrt{n_x} L_2 \tau + \frac{2\sqrt{n_x}\Delta}{\tau} \right)\nonumber\\ &&+ 9\gamma\left(\sigma^2 + n_x L_2^2 \tau^2 + \frac{2n_x\Delta^2}{\tau^2}\right).
\end{eqnarray}
\end{theorem}
\textbf{Proof of }\eqref{teor5.1}:
We begin with \eqref{tlemma} and taking $z = z_k$, $g = \gamma\widetilde{g}_f(z_k, e_k, \tau, \xi_k)$, $g_{1/2} = \gamma \widetilde{g}_f(z_{k+1/2}, e_{k+1/2}, \tau, \xi_{k+1/2})$, then $z_{1/2} = z_{k+1/2}$, $z_{1} = z_{k+1}$ and we get

\begin{eqnarray*}
    \gamma \langle \widetilde{g}_f(z_{k+1/2}, e_{k+1/2}, \tau, &\xi_{k+1/2}&), z_{k+1/2} - u \rangle  \\
     &\leq&  V_{z_k}(u) - V_{z_{k+1}}(u) -  V_{z_k}(z_{k+1/2}) \nonumber \\
     &&+ \frac{\gamma^2}{2}\|\widetilde{g}_f(z_{k+1/2}, e_{k+1/2}, \tau, \xi_{k+1/2}) - \widetilde{g}_f(z_k, e_k, \tau, \xi_k)\|^2_q \nonumber\\
      &\overset{\eqref{eq:squared_sum}}{\leq}& V_{z_k}(u) - V_{z_{k+1}}(u) -  V_{z_k}(z_{k+1/2}) \nonumber \\ 
      &&+ \frac{3\gamma^2}{2}\|F(z_{k+1/2}) - F(z_k)\|^2_q \nonumber \\ 
      &&+ \frac{3\gamma^2}{2}\|\widetilde{g}_f(z_{k+1/2}, e_{k+1/2}, \tau, \xi_{k+1/2}) - F(z_{k+1/2})\|^2_q \nonumber \\ 
     &&+ \frac{3\gamma^2}{2}\|\widetilde{g}_f(z_k, e_k, \tau, \xi_k) - F(z_k)\|^2_q 
\end{eqnarray*}  
With \eqref{SP_c} it gives
\begin{eqnarray*}
\gamma \langle \widetilde{g}_f(z_{k+1/2}, e_{k+1/2}, \tau, &\xi_{k+1/2}&), z_{k+1/2} - u \rangle  \\
    &{\leq}& V_{z_k}(u) - V_{z_{k+1}}(u) -  V_{z_k}(z_{k+1/2}) \nonumber \\ 
     &&+ \frac{3\gamma^2L^2}{2}\|z_{k+1/2} - z_k\|^2_2 \nonumber \\ 
     &&+ \frac{3\gamma^2}{2}\|\widetilde{g}_f(z_{k+1/2}, e_{k+1/2}, \tau, \xi_{k+1/2}) - F(z_{k+1/2})\|^2_q \nonumber \\ 
     &&+ \frac{3\gamma^2}{2}\|\widetilde{g}_f(z_k, e_k, \tau, \xi_k) - F(z_k)\|^2_q.
\end{eqnarray*}
Applying the property: $V_{z_{k}}(z_{k+1/2}) \geq \nicefrac{1}{2} \|z_{k+1/2} - z_k \|^2 \geq \nicefrac{1}{2} \|z_{k+1/2} - z_k \|^2_2$, with $\gamma \leq \nicefrac{1}{2L}$, we get
\begin{eqnarray*}
    \gamma \langle \widetilde{g}_f(z_{k+1/2}, e_{k+1/2}, &\tau,& \xi_{k+1/2}), z_{k+1/2} - u \rangle  \leq  V_{z_k}(u) - V_{z_{k+1}}(u) \nonumber \\ 
    &&+ \frac{3\gamma^2}{2}\|\widetilde{g}_f(z_{k+1/2}, e_{k+1/2}, \tau, \xi_{k+1/2}) - F(z_{k+1/2})\|^2_q \nonumber \\ 
     &&+ \frac{3\gamma^2}{2}\|\widetilde{g}_f(z_k, e_k, \tau, \xi_k) - F(z_k)\|^2_q. 
\end{eqnarray*}
Taking the full expectation and using \eqref{lemma1_3_a}, \eqref{lemma1_3_1_a} with \eqref{b}:
\begin{eqnarray*}
    \EE\left[\gamma \langle  F(z_{k+1/2}), z_{k+1/2} - u \rangle \right] &\leq&  \EE\left[V_{z_k}(u) \right]- \EE\left[V_{z_{k+1}}(u)\right] \nonumber \\ 
    && + 2\gamma  \left(\sqrt{n_x} L_2 \tau + \frac{2\sqrt{n_x}\Delta}{\tau} \right)D_p\nonumber\\
    &&+ 3\gamma^2\left(3\sigma^2 + 3n_x L_2^2 \tau^2 + \frac{6n_x\Delta^2}{\tau^2}\right).
\end{eqnarray*}
It remains to sum up from $k=0$ to $k=N$ and use \ref{connect_e_sad} and finish the proof of this theorem.
\EndProof

\end{document}